\documentclass{amsart}

\usepackage{adjustbox}
\usepackage[colorlinks=true, allcolors=blue]{hyperref}
\usepackage{amsmath}
\usepackage{amstext} 
\usepackage{array}   
\usepackage{mathtools}
\usepackage{tikz}
\usepackage{tikz-cd}

\usepackage{amsmath}
\usepackage{amssymb}
\usepackage{amsfonts}
\usepackage{amsthm}

\theoremstyle{plain}
\newtheorem*{thm*}{Theorem}
\newtheorem{thm}{Theorem}[section]
\newtheorem*{lemma*}{Lemma}

\newtheorem*{pro*}{Proposition}

\newtheorem*{cor*}{Corollary}
\newtheorem*{con*}{Conjecture}
\newtheorem{con}[thm]{Conjecture}
\newtheorem*{remark*}{Remark}
\newtheorem{remark}[thm]{Remark}

\theoremstyle{definition}

\newtheorem*{df*}{Definition}

\theoremstyle{remark}
\newtheorem{rem}[thm]{Remark}
\newtheorem*{rem*}{Remark}
\newtheorem{ex}[thm]{Example}
\newtheorem*{ex*}{Example}

\newcolumntype{L}{>{$}l<{$}} 
\newcolumntype{C}{>{$}c<{$}} 

\def\SK{\hskip 1 true cm }
\def\HOT{h.o.t.}
\def\C{\mathbf C}
\def\R{\mathbf R}
\def\Q{\mathbf Q}
\def\PP{\mathbf P}
\def\csm{c^{sm}}
\def\ssm{s^{sm}}
\def\HH{H^{**}}
\def\HD{H^{\leq 2d}}
\def\ssmD{\ssm_{\leq d}}
\DeclareMathOperator{\codim}{codim}

\def\E{\mathcal E}
\def\K{\mathcal K}
\def\A{\mathcal A}
\def\O{\mathcal O}
\def\Ideal{\mathcal I}
\DeclareMathOperator\Diff{Diff}
\DeclareMathOperator{\Thom}{Thom}
\DeclareMathOperator{\GL}{GL}
\DeclareMathOperator{\T}{T}
\def\rhosou{\rho^{\text{source}}}
\def\rhotar{\rho^{\text{target}}}

\title[Interpolation characterization of higher Thom polynomials]{Interpolation characterization of \\ higher Thom polynomials}
\author{Rich\'ard Rim\'anyi}

\begin{document}

\begin{abstract}
Thom polynomials provide universal formulas for the fundamental class of singularity loci in terms of characteristic classes. Ohmoto extended this notion to SSM-Thom polynomials, which refine this description by capturing the richer Segre-Schwartz-MacPherson (SSM) class of singularity loci. While previous methods for computing SSM-Thom polynomials relied on intricate geometric arguments, we introduce a more efficient approach that depends solely on the symmetries of singularities. Our method is inspired by connections to Geometric Representation Theory, particularly the interpolation properties of Maulik-Okounkov stable envelopes. By formulating SSM analogs of these axioms within a degree-bounded framework, we obtain new computational tools for SSM-Thom polynomials. We also present explicit examples of SSM-Thom polynomials, and illustrate their applications in enumerative geometry and singularity theory.
\end{abstract}

\maketitle

\section{Introduction}

Thom polynomials are universal formulas that describe geometric information of singularity loci in terms of simple invariants of a map.

Specifically, consider a complex analytic map $F:M^m \to N^n$ between compact complex manifolds. Given a singularity type $\eta$, the singularity locus $\eta(F)$ is the subset of $M$ where $F$ has singularity $\eta$. A key result is that there exists a multivariate polynomial, known as the Thom polynomial, which depends only on $\eta$ and satisfies the following property: for any (non-degenerate) map $F$, the fundamental class of 
$\overline{\eta(F)}$ is given by evaluating the Thom polynomial at the characteristic classes of the map \cite{thomoriginal, primer}.

Characteristic classes are generally easy to compute, making Thom polynomials a valuable tool for extracting refined information about degeneracy loci, such as their number or (multi-)degree. Today, powerful techniques exist for computing Thom polynomials for a wide range of singularities, drawing on methods from geometry or algebra, including resolution of singularities, embedded resolution, Gr\"obner degeneration, Hecke-type recursions, iterated residues, and interpolation.

T. Ohmoto’s pioneering works \cite{OhmotoCamb,ohmotoSMTP} extended the notion of Thom polynomials to {\em higher Thom polynomials}, also known as {\em SSM-Thom polynomials}. An SSM-Thom polynomial is an inhomogeneous power series in the same variables as the Thom polynomial, with its lowest-degree term being the Thom polynomial itself. Beyond capturing the fundamental class, the full SSM-Thom polynomial encodes the Segre-Schwartz-MacPherson (SSM) class of a singularity locus.

SSM classes provide significantly more refined information about subvarieties in a smooth ambient space than fundamental classes alone. For instance, in the case of a curve in $\PP^3$, its fundamental class encodes its degree, while its SSM class encodes both its degree and its Euler characteristic.

Beyond establishing the definition and key structural theorems, Ohmoto et al. \cite{ohmotoSMTP, nekarda1, nekarda2}, also developed a method to compute some initial terms of certain SSM-Thom polynomials. That approach relies on the challenging geometric task of identifying the $\eta$-singularity loci of so-called prototypes of other singularities.

The aim of this paper is to introduce a more effective method that bypasses these geometric considerations, relying instead solely on the determination of symmetries of singularities. This approach is formulated in our main Theorem~\ref{thm:main}.

The key idea behind our main theorem is a useful connection to Geometric Representation Theory. In that field, Maulik and Okounkov \cite{MO} introduced the concept of {\em stable envelope}, which have since been successfully applied in various contexts. Notably, in certain cases and under appropriate identifications, stable envelopes can be interpreted as Chern-Schwartz-MacPherson (CSM) classes of specific subvarieties \cite{RV, AMSS_csm}. Moreover, the defining interpolation properties of stable envelopes can be adapted to provide an interpolation-based characterization of CSM classes for orbits of certain representations \cite{FRcsm}. Our proof follows this framework by formulating and employing SSM versions of interpolation axioms constrained by a degree bound.

After the proof (Sections~\ref{sec:proof1} and \ref{sec:proof2}) in Sections~\ref{sec:expansions} and \ref{sec:app} we study the algebraic combinatorics, and present three concrete examples showcasing applications of the newly computed SSM-Thom polynomials. Throughout the paper, we give many examples, and collect many others on the Thom Polynomial Portal \cite{TPP}.

\bigskip

\noindent{\bf Acknowledgment.} The author was supported by the National Science Foundation under Grant No. 2200867. Special thanks to L.~Feh\'er and J. Koncki for helpful discussions during the preparation of this paper.

\section{Chern-Schwartz-MacPherson theory} \label{sec:CSM}

Let $M$ be a smooth, compact complex ambient variety, and let $\Sigma\subset M$ be a possibly singular subvariety of codimension $k$. The concept of the fundamental class $[\Sigma \subset M] \in H^{2k}(M)$\footnote{In the whole paper cohomology is meant with rational coefficients} admits a deformation which comes in two distinct versions: the Chern-\-Schwartz-\-Mac\-Pher\-son (CSM) class and the Segre-Schwartz-\-Mac\-Pher\-son (SSM) class. These two versions differ only by an explicit multiplicative factor:
\[
\begin{array}{rcl}
\csm(\Sigma\subset M)&=[\Sigma \subset M]+\HOT & \in H^*(M),
\\
\ssm(\Sigma\subset M)=\frac{\csm(\Sigma\subset M)}{c(TM)} & =[\Sigma \subset M] + \HOT & \in H^{*}(M).
\end{array}
\]
When the ambient variety is clear, we just write $\csm(\Sigma)$, $\ssm(\Sigma)$, $[\Sigma]$. If a group $G$ acts on $M$ and $\Sigma$ is an invariant subvariety, the CSM and SSM classes exist in equivariant cohomology. The equivariant SSM class may have arbitrary large degree components, hence it is an element of the completion $H_G^{**}(M)$ of $H_G^*(M)$.

For us the properties (i)-(iv) below will serve as the definition of CSM and SSM classes. In fact, these properties overdetermine the concepts. Historically, the existence of CSM classes satisfying conditions (i), (iii) and (iv) was conjectured by Grothendieck---as a substitute for {\em total Chern class of the tangent bundle} of $\Sigma$ when $\Sigma$ is singular. This conjecture was proved independently by Schwartz and MacPherson \cite{BrSch,macpherson}. The equivariant version is due to Ohmoto \cite{OhmotoCamb}. 

\begin{enumerate}
\item[(i)] The CSM class is consistent with push-forward: If $f : Y \to M$ is an equivariant morphism between smooth varieties then $f_*(c(TY)) = \sum_a \csm(M_a)$, where $M_a = \{m \in M : \chi(f^{-1}(m)) = a\}$.
\item[(ii)] The SSM class is consistent with pull-back: For $Y$ and $M$ smooth let $X\subset M$ be an invariant subvariety with an invariant Whitney stratification. Assume $f: Y \to M$ is transversal to the strata of $X$. Then $\ssm(f^{-1}(X)) = f^*(\ssm(X))$.
\item[(iii)] CSM classes satisfy normalization: if $i : \Sigma \subset M$ is a smooth closed subvariety then
$\csm(\Sigma \subset M ) = i_*(c(T \Sigma))$.
\item[(iv)] Both CSM and SSM classes are defined for all constructible functions on the ambient space, where
the $\csm(\Sigma \subset M )$ concept corresponds to applying it to the indicator function of $\Sigma$.
In this generality CSM classes (hence SSM classes too) are additive:
$\csm(\alpha f + \beta g) = \alpha \csm(f ) + \beta \csm(g)$. In fact, a precise definition of the CSM class should start with the additivity property, and in this language CSM is a natural transformation between the functor of constructible functions and (Borel-Moore) homology. Our version is its Poincar\'e dual in the ambient space.
\end{enumerate}

What geometric information does the CSM class of a subvariety encode? To illustrate the answer, assume that the ambient space is a projective space $\PP^N$. In that case the lowest degree part, the fundamental class, is the degree of $\Sigma\subset \PP^N$. The higher degree components encode (in an algebraic way that we will recall precisely in Section~\ref{sec:application}) the Euler characteristics of the various dimensional linear sections of $\Sigma$. For example, if $\Sigma$ is a connected curve in $\PP^N$ then the two terms of its CSM class encode its degree and its genus.

\section{Contact singularities and their global theory}
In this section we give a brief review on the local and global theory of contact singularities of maps. We recommend the recent book \cite{MB} and the surveys \cite{OhmotoSurvey,primer} as references. 

\subsection{Contact singularities} \label{sec:contact}
Let $m,\ell\geq 0$. Let $\E(m,m+\ell)$ be the vector space of germs of holomorhic maps $(\C^m,0)\to (\C^{m+\ell},0)$. 

The group of complex holomorphic diffeomorphism germs of $(\C^m\times \C^{m+\ell},0)$ of the form 
\[
\Phi(x,y)=(\phi(x), \psi(x,y))
\]
where $\psi(x,0)=0$, is called the contact group $\K(m,m+\ell)$. The left-right group $\A(m,m+\ell)=\Diff(\C^m,0)\times \Diff(\C^{m+\ell},0)\subset \K(m,m+\ell)$ is the product of the groups of holomorphic diffeomorphism germs of the source and the target.

The group $\K(m,m+\ell)$ acts on $\E(m,m+\ell)$ via its action on the graph. Orbits of the action are called contact singularities, orbits of $\A(m,m+\ell)$ are called right-left singularities.

The local algebra of a germ $f:(x_1,\ldots,x_m)\mapsto (f_1,\ldots,f_{m+\ell})$ is defined as 
$Q_f=\O_m/f^*\Ideal_{m+\ell}$, where $\O_m$ is the ring of holomorphic function germs at $(\C^m,0)$ and $\Ideal_m$ is its maximal ideal.
We will be only interested in finite germs, that is, when this algebra is finite dimensional and can also be presented as $Q_f=\C[[x_1,\ldots,x_m]]/(f_1,\ldots,f_{m+\ell})$. It is a theorem of Mather \cite{mather4} that two germs in $\E(m,m+\ell)$ are $\K(m,m+\ell)$-equivalent if and only if their local algebras are isomorphic. Hence a (commutative, finite dimensional, local) algebra $Q$ as well as $m$ and $\ell$ determine an orbit $\eta(Q,m,\ell)$ (unless this set is empty).

\smallskip

\noindent{\em Notation.} As is standard in singularity theory, to save space, we will often use the ideal $J$ itself to denote the algebra $\C[x_i]/J$. For example, we will write $(x^2,y^3)$, instead of the algebra: $\C[x,y]/(x^2,y^3)$. In these cases the omitted ring is always a polynomial ring in the variables that turn up in the ideal.


\begin{rem}
    In practice, we can replace the vector space $\E(m,m+\ell)$ of germs, and the group $\K(m,m+\ell)$ of germs with their $N$-jets ($N \gg 0$) to obtain $\E_N(m,m+\ell)$ and $\K_N(m,m+\ell)$. This way we obtain an algebraic group acting on a finite dimensional vector space. Our constructions and results do not depend on $N$ as long as $N$ is large enough, hence by abuse of notation we will not write the subscript~$N$.
\end{rem}

Define the {\em Mather-bound}
    \[
    M(\ell)=\begin{cases} 6\ell+8 & \text{if } \ell=0,1,2,3 \\ 6\ell+7 & \text{if } \ell\geq 4. \end{cases}
    \]
Dimension pairs $(m,m+\ell)$ satisfying $m\leq M(\ell)$ are called the `nice dimensions', see \cite{mather6, MB}. The celebrated result of Mather (not explicitly used in this paper) is that so-called {\em stable} singularities are dense exactly in the nice dimensions. 

\begin{thm}
   
    \begin{itemize}
    \item For every $m$ there exists a $\K(m,m+\ell)$ invariant subvariety $\Pi(m,m+\ell) \subset \E(m,m+\ell)$ of codimension $>M(\ell)$, such that $\E(m,m+\ell)$ has finitely many $\K(m,m+\ell)$ orbits, all of codimension $\leq M(\ell)$.
    \item The list of algebras corresponding to the finitely many orbits in $\E(m,m+\ell)-\Pi$ is an explicit list given in \cite{mather6, MB}.
    \item For a fixed $\ell$ this list does not depend on $m$ for large enough $m$.
    \end{itemize}
\end{thm}

In essence, the theorem claims that for a given $\ell$ the classification of $\K(*,*+\ell)$ orbits of $\E(*,*+\ell)$ is finite and is explicitly known, up to codimension $M(\ell)$. A further remarkable property of the finitely many singularities occurring in codimension $\leq M(\ell)$ will be discussed in Section~\ref{sec:proof2}.


\begin{rem} \label{rem:ListOfAlgebras}
We claimed that the list of algebras is explicitly provided in the literature. While this statement is essentially accurate, it comes with some caveats. First, the nature of such huge lists is that some typos occur---even if they are easy to spot for the experts. Moreover, the cited literature lists the {\em real forms} of the algebras, hence for us the list contains repetitions. Some of those repetitions are obvious, for example over $\C$
    \[ 
    (xy,x^2+y^2) \cong (xy,x^2-y^2) \cong (x^2,y^2), \text{ or}
    \]
    \[
    (x^2+y^2,x^3) \cong (xy,x^3+y^3),
    \]
    some others are less obvious, like
    \[
    (x^2+y^2+z^2,xy,xz,yz) \cong (x^2,y^2,z^2-xy,xz+yz).
    \]
    Remarkably, the global theory of singularities provides a powerful tool for addressing such questions; see Section~\ref{sec:RC}.
\end{rem}

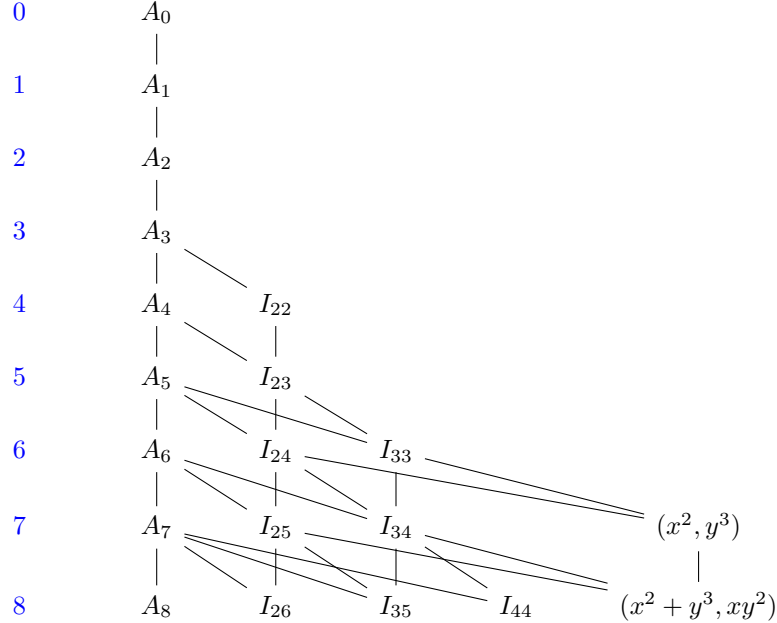
\begin{figure}
\[
\begin{tikzcd}[column sep=24pt, row sep=12pt]
{\color{blue} codim} & \\
{\color{blue}0} & A_0\arrow[-,d] \\
{\color{blue}1} & A_1\arrow[-,d] \\
{\color{blue}2} & A_2\arrow[-,d] \\
{\color{blue}3} & A_3\arrow[-,d]\arrow[-,dr] \\
{\color{blue}4} & A_4\arrow[-,d]\arrow[-,dr] & I_{22}\arrow[-,d]\\
{\color{blue}5} & A_5\arrow[-,d]\arrow[-,dr]\arrow[-,drr] & I_{23}\arrow[-,d]\arrow[-,dr]\\
{\color{blue}6} & A_6\arrow[-,d]\arrow[-,dr]\arrow[-,drr] & I_{24}\arrow[-,d]\arrow[-,dr]\arrow[-,drrr] & I_{33}\arrow[-,d]\arrow[-,drr]\\
{\color{blue}7} & A_7\arrow[-,d]\arrow[-,dr]\arrow[-,drr]\arrow[-,drrr] & I_{25}\arrow[-,d]\arrow[-,dr]\arrow[-,drrr] & I_{34}\arrow[-,d]\arrow[-,dr]\arrow[-,drr] & & (x^2,y^3)\arrow[-,d]\\
{\color{blue}8} & A_8                        & I_{26}                        & I_{35} & I_{44} & (x^2+y^3,xy^2)  \\
\end{tikzcd}
\]
\caption{Classification of contact singularities for $\ell=0$ in the nice dimensions. Edges indicate the covering relation of the adjacency hierarcy.}
\label{fig:l0}
\end{figure}

\begin{ex}
    For $\ell=0$ the Mather-bound is $8$. The complete and repetition-free list of algebras corresponding to the codim$\leq 8$ orbits in $\E(m,m+\ell)$ ($m\gg 0$) is
    \begin{itemize}
        \item $A_i=\C[x]/(x^{i+1})$ for $i=0,\ldots,8$,
        \item $I_{a,b}=\C[x,y]/(xy,x^a+y^b)$ for $2\leq a\leq b$, $a+b\leq 8$,
        \item $\C[x,y]/(x^2,y^3)$,
        \item $\C[x,y]/(x^2+y^3,xy^2)$.
    \end{itemize}
    The codimensions and the adjacency-hierarchy  of the orbits $\eta(Q,m,0)$ for $m\gg 0$ (namely, which orbits are contained in which other orbits' closures) are displayed in Figure~\ref{fig:l0}.
\end{ex}

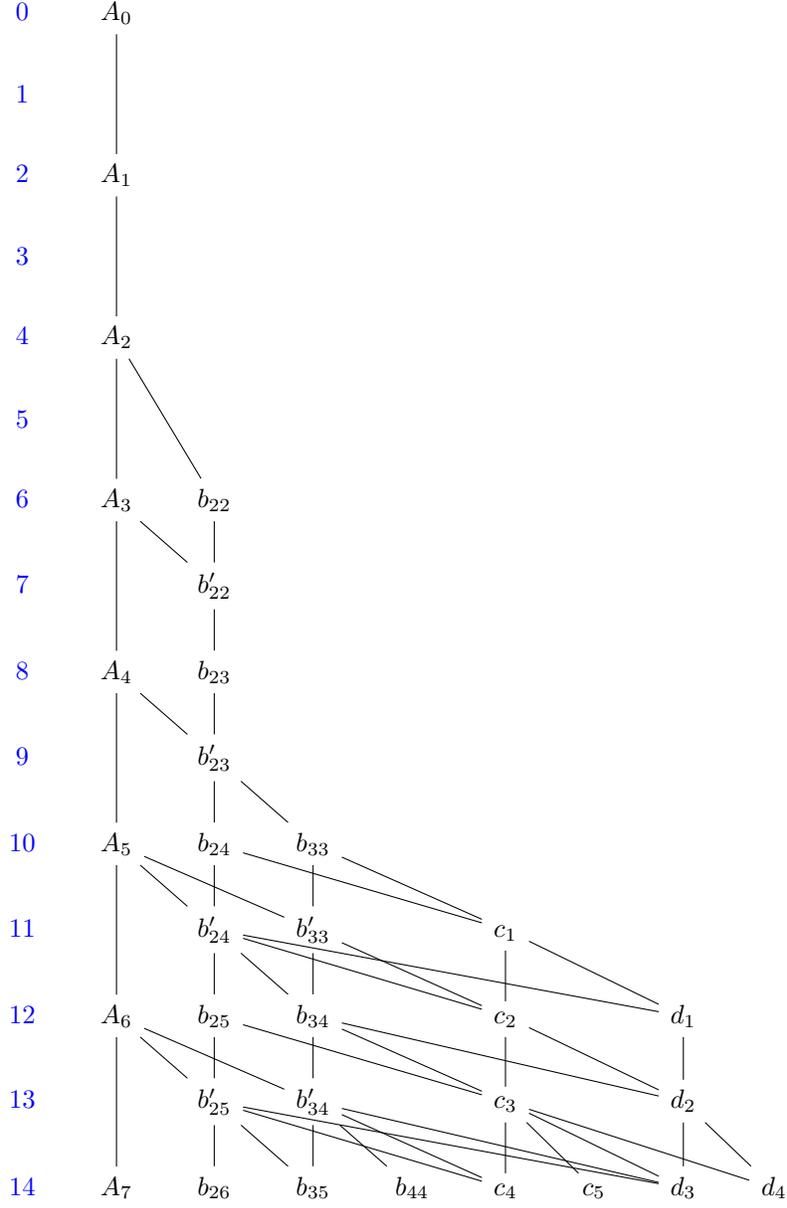
\begin{figure}
\[
\begin{tikzcd}[column sep=16pt, row sep=16pt]
{\color{blue}0} &A_0\ar[-,dd] \\
{\color{blue}1} & \\
{\color{blue} 2} &A_1\ar[-,dd] \\
{\color{blue} 3} & \\
{\color{blue}4} &A_2\ar[-,dd]\ar[-,ddr] \\
{\color{blue}5} &\\
{\color{blue}6} &A_3\ar[-,dd]\ar[-,dr] & b_{22}\ar[-,d] \\
{\color{blue}7} &    & b'_{22}\ar[-,d] \\
{\color{blue}8} &A_4\ar[-,dd]\ar[-,dr] & b_{23}\ar[-,d] \\
{\color{blue}9} &    & b'_{23}\ar[-,d]\ar[-,dr] \\
{\color{blue}10} &A_5\ar[-,dd]\ar[-,dr]\ar[-,drr] & b_{24}\ar[-,d]\ar[-,drrr] & b_{33}\ar[-,d]\ar[-,drr] \\
{\color{blue}11} &    & b'_{24}\ar[-,d]\ar[-,dr]\ar[-,drrr]\ar[-,drrrrr] &b'_{33}\ar[-,d]\ar[-,drr] & & c_1\ar[-,d]\ar[-,drr] \\
{\color{blue}12} &A_6\ar[-,dd]\ar[-,dr]\ar[-,drr] & b_{25}\ar[-,d]\ar[-,drrr] & b_{34}\ar[-,drr]\ar[-,drrrr]\ar[-,d] & & c_2\ar[-,d]\ar[-,drr] & & d_1\ar[-,d] \\
{\color{blue}13} &    & b'_{25}\ar[-,d]\ar[-,dr]\ar[-,drrr]\ar[-,drrrrr] & b'_{34}\ar[-,d]\ar[-,dr]\ar[-,drr]\ar[-,drrrr] & & c_3\ar[-,d]\ar[-,dr]\ar[-,drr]\ar[-,drrr] & & d_2\ar[-,d]\ar[-,dr]\\
{\color{blue}14} &A_7 & b_{26} & b_{35} & b_{44} & c_4 & c_5 & d_3 & d_4
\end{tikzcd}
\]
\caption{Classification and hierarchy of contact singularities for $\ell=1$ in the nice dimensions. Notation: $b_{ij}=III_{ij}$, $b'_{ij}=I_{ij}$, $c_1=(x^2,xy^2,y^3)$, $c_2=(x^2,y^3)$, $c_3=(x^2+y^3,xy^2,y^4)$, $c_4=(x^2+y^3,xy^2)$, $c_5=(x^2,xy^2,y^4)$, $d_1=(x^2+y^2+z^2,xy,xz,yz)$, $d_2=(x^2,y^2,z^2,xy+xz)$, $d_3=(x^2-y^2+z^3,xy,xz,yz)$, $d_4=(x^2+yz,xz,y^2,z^2)$.}
\label{fig:l1}
\end{figure}

\begin{ex}
    For $\ell=1$ we have $M(\ell)=14$. The complete and repetition-free list of algebras corresponding to the codim$\leq 14$ orbits in $\E(m,m+\ell)$ ($m\gg 0$), together with their hierarchy is given in Figure~\ref{fig:l1}.
\end{ex}

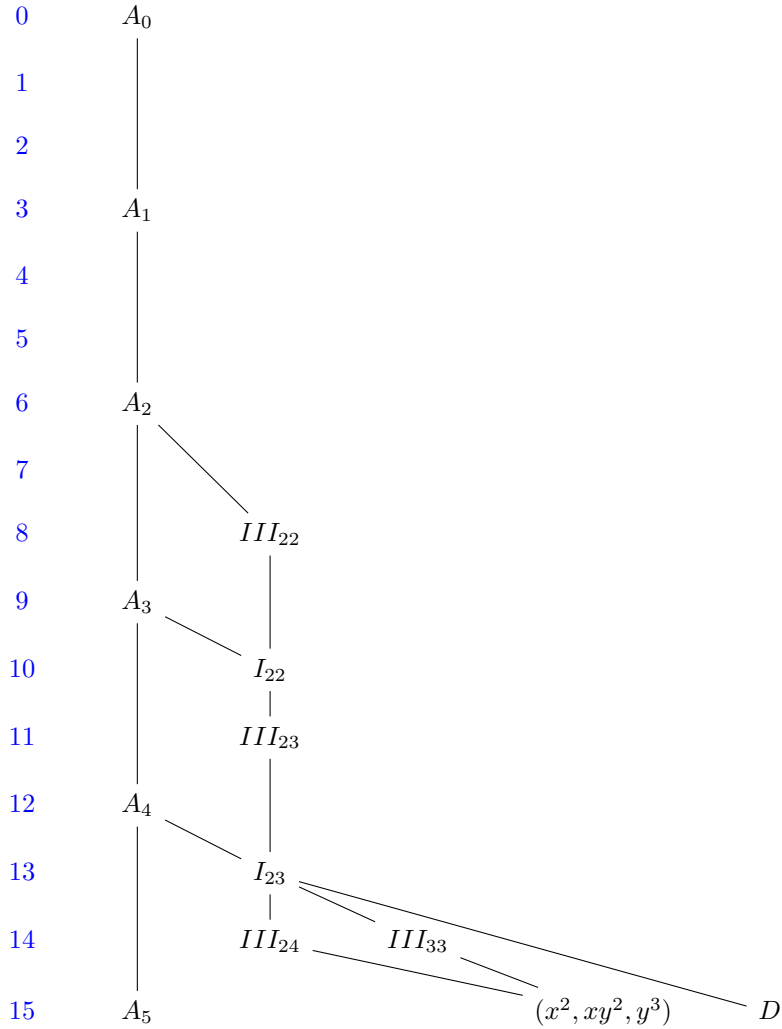
\begin{figure}  
\[
\begin{tikzcd}[column sep=24pt,row sep=10pt]
{\color{blue}0} & A_0\ar[-,ddd] &  \\
{\color{blue}1} & & & \\
{\color{blue}2} & & &  \\
{\color{blue}3} & A_1\ar[-,ddd] & & \\
{\color{blue}4} & & &   \\
{\color{blue}5} & & & &  \\
{\color{blue}6} & A_2\ar[-,ddd]\ar[-,ddr]& &  \\
{\color{blue}7} & & &\\
{\color{blue}8} & & III_{22}\ar[-,dd] &  \\
{\color{blue}9} & A_3\ar[-,ddd]\ar[-,dr] & & \\
{\color{blue}10} & & I_{22}\ar[-,d]  & \\
{\color{blue}11} & & III_{23}\ar[-,dd] & \\
{\color{blue}12} & A_4\ar[-,ddd]\ar[-,dr] & & \\
{\color{blue}13} & & I_{23}\ar[-,d]\ar[-,dr]\ar[-,ddrrr] & \\
{\color{blue}14} & & III_{24}\ar[-,drr] & III_{33}\ar[-,dr]\\
{\color{blue}15} & A_5 & & & (x^2,xy^2,y^3) & D
\end{tikzcd}
\]
\caption{Classification of contact singularities for $\ell=2$ up to codimension 15. Here $D$ stands for $(x^2-y^2,x^2-z^2,xy,xz,yz)$. The classification is finite up to codimension $M(3)=20$. The classification of those with codimension$\in[16,20]$ are not shown.}
\label{fig:l2}
\end{figure}

\begin{ex}
    For $\ell=2$ the Mather-bound is $20$. Figure~\ref{fig:l2} lists the algebras corresponding to orbits up to codimension $15$ only. 
\end{ex}
    
\subsection{Genotypes, prototypes, and symmetries of singularities}\label{sec:geno_etc}

Let $Q$ be a finite dimensional, commutative, local algebra; and let 
\[
Q=\C[x_1,\ldots,x_a]/(r_1,\ldots,r_b)
\]
be a presentation with minimal $b-a$ and $a$. For $\ell\geq b-a$ the germ
\[
(\C^a,0)\to (\C^{a+\ell},0),
\qquad
(x_1,\ldots,x_a) \mapsto (r_1,\ldots,r_b,\underbrace{0,\ldots,0}_{a+\ell-b})
\]
is called the {\em genotype} for the singularity with local algebra $Q$ and relative dimension~$\ell$. A procedure called `(mini)versal unfolding' \cite[\S5]{MB}, \cite[\S1.1]{rrtp} applied to the genotype produces what is called the {\em prototype} $p:(\C^s,0)\to (\C^{s+\ell},0)$ for the singularity with local algebra $Q$ and relative dimension $\ell$. 

\begin{ex} \label{ex:x2y3}
    Let $Q=\C[x,y]/(x^2,y^3)$. For $\ell=0$ and $\ell=1$ the genotypes are 
    \[ (x,y)\mapsto (x^2,y^3),
    \qquad
    (x,y)\mapsto (x^2,y^3,0).
    \]
    The respective prototypes are
    \begin{multline*}
    (\C^{7},0)\to (\C^{7},0),\\
      (x,y,u_1,\ldots,u_{5})\mapsto 
        (x^2+u_1y+u_2y^2, y^3+u_3x+u_4y+u_5xy,u_1,\ldots,u_{5}),
    \end{multline*}
    \begin{multline*}   
   (\C^{12},0)\to (\C^{13},0)\\
       (x,y,u_1,\ldots,u_{10})\mapsto 
        (x^2+u_1y+u_2y^2, y^3+u_3x+u_4y+u_5xy, \\ 
        u_6x+u_7y+u_8xy+u_9y^2+u_{10}xy^2,u_1,\ldots,u_{10}).
    \end{multline*}
\end{ex}
The maximal compact (MC) subgroup of the the symmetry group of $p$
\[
G=MC\{ (g_1,g_2)\in \Diff(\C^s,0) \times \Diff(\C^{s+\ell},0) :
 g_2 \circ p \circ g_1^{-1} =p \},
\]
will play a central role in our main Theorem \ref{thm:main}. Specifically, the genotype and prototype can be chosen such that $G\subset \GL_s(\C)\times \GL_{s+\ell}(\C)$. Furthermore, the only role of $G$ throughout the paper is considering the kernels of certain ring homomorphisms into $H^*(BG)$. Since $H^*(BG)\subset H^*(BT)$ (by the splitting lemma) for the maximal torus $T\leq G$, we can equivalently assume that $G$ denotes not the symmetry group but its maximal torus. 

The linear representations of $G$ on $\C^{s}$ and $\C^{s+\ell}$ will be called $\rhosou$ and $\rhotar$. We use $e(\rho)$ and $c(\rho)$ to denote the Euler and total Chern classes of a representation~$\rho$.
The formal power series
\[
c(Q,\ell)=\frac{ c(\rhotar) }{c(\rhosou)} \in H^{**}(BG).
\]
will play a role later.

\begin{ex}
    Continuing Example~\ref{ex:x2y3}, for $\ell=0$  we have
    \[
    G=U(1)\times U(1), \qquad
    \rho_0=\lambda_1 \oplus \lambda_2 \oplus \rho_U, \qquad
    \rho_1=\lambda_1^2 \oplus \lambda_2^3 \oplus \rho_U,
    \]
    where $\lambda_i$ is the standard representation of the $i$'th factor of $G$, and 
    \[
    \rho_U= \lambda_1^2\lambda_2^{-1} \oplus \lambda_1^2\lambda_2^{-2} \oplus \lambda_2^3\lambda_1^{-1} \oplus \lambda_2^2 \oplus \lambda_2^2\lambda_1^{-1}.  
    \]
    Using the notation $H^*(BG)=\Q[a,b]$ we have
    \begin{multline*}
    c(Q,0)=\frac{c(\rhotar)}{c(\rhosou)}=\frac{c(\lambda_1^2 \oplus \lambda_2^3)c(\rho_U)}{c(\lambda_1 \oplus \lambda_2)c(\rho_U)}=
    \frac{(1+2a)(1+3b)}{(1+a)(1+b)}=\\
    1+(a+2b)+(-a^2+2ab-2b^2)+(a^3-2a^2b-2ab^2+2b^3)+\ldots
    \end{multline*}
    
\end{ex}

\subsection{The global theory}\label{sec:global}
Let $\eta$ be a $\K(m,m+\ell)$ invariant set in $\E(m,m+\ell)$.  Since $\A(m,m+\ell)\subset \K(m,m+\ell)$, the definition of the singularity locus 
\[
\eta(F)=\{ x\in M^m : \text{the germ of $F$ at $x$ is in $\eta$} \}
\]
makes sense for a map $M^m\to N^{m+\ell}$ between complex manifolds. 

\subsubsection{Universal formulas for degeneracy loci}
The global theory of singularities seeks `universal formulas' for various characteristics of the singularity locus $\eta(F)$, in terms of simple invariants of the map $F:M\to N$. A key example for such universal formulas is the Thom polynomial. 

\begin{thm} \label{thm:ThomPoly}\cite{thomoriginal,OhmotoSurvey}
    Let $\eta$ be a $\K(m,m+\ell)$ invariant subset of $\E(m,m+\ell)$ of pure codimension $k$. There exists a polynomial of degree $k$, the Thom polynomial of $\eta$, 
    \[
    \Thom(\eta)\in \Q[c_1,c_2,\ldots] \qquad\qquad (\deg(c_i)=i)
    \]
    such that for all maps $F:M^m\to N^{m+\ell}$ satisfying a transversality condition we have
    \begin{equation} \label{eq:Tp=locus}
     [\overline{\eta(F)}] = \Thom(\eta)( c(F) ).
     \end{equation}
    Moreover, the Thom polynomial of $\eta(Q,m,\ell)$ depends only on $Q$ and $\ell$, not on $m$. This Thom polynomial will be denoted by $\Thom(Q,\ell)$.
\end{thm}

In addition to the original paper \cite{thomoriginal} and the recent survey article \cite{OhmotoSurvey}, the reader may find \cite{damon,rrtp,kaza:noas,FRannals,BSz,BercziNew} instructive.

Recall that $[\Sigma]=[\Sigma\subset M]\in H^{2k}(M)$ is the cohomological fundamental class represented by $\Sigma\subset M$. The Chern classes of a map $F:M\to N$ are defined by
\[
c(F)=1+c_1(F)+c_2(F)+\ldots = \frac{F^*(c(TN))}{c(TM)}= \frac{c(F^*(TN))}{c(TM)}.
\]

The power of the theorem is that $\Thom(\eta)$ does not depend on the concrete situation, that is, on $M$, $N$ or $F$. It is `universal': it only depends on the singularity~$\eta$. The concrete situation of $F:M\to N$ only enters \eqref{eq:Tp=locus} via the simple characteristic classes of $F$.

\subsubsection{Equivariant fundamental class interpretation}
\label{sec:TPequivariantclass}
    Let $G_{m,m+\ell}=\GL_m(\C)\times \GL_{m+\ell}(\C)$ with its conjugation action on $\E(m,m+\ell)$. Reference to the transversality conditions in Theorem~\ref{thm:ThomPoly} can be avoided if we interpret $\Thom({\eta})$ as the $G_{m,m+\ell}$-equivariant fundamental class of $\overline{\eta} \subset \E(m,m+\ell)$: 
\[
\Thom(\eta)=[ \overline{\eta}] \in H_{G_{m,m+\ell}}^*(\E(m,m+\ell))= H_{G_{m,m+\ell}}^*(pt).
\]
In this interpretation the variables $c_i$ are defined as
\begin{equation}
\label{eq:quotient_abstract}
1+c_1+c_2+\ldots=\frac{1+b_1+b_2+\ldots+b_{m+\ell}}{1+a_1+a_2+\ldots+a_m},
\end{equation}
where $a_i$ resp. $b_i$ are the Chern classes associated to the factors $\GL_m(\C)$ resp. $\GL_{m+\ell}(\C)$ of $G_{m,m+\ell}$. In this interpretation the `universality' property mentioned above is a tautology.

\smallskip

The theory of Thom polynomials has been intensively studied and applied in geometry and topology; see the recent surveys \cite{primer, OhmotoSurvey}. In this paper we study a 1-parameter deformation of Thom polynomials, the SSM-Thom polynomials.

\section{SSM-Thom polynomials by interpolation}

\subsection{SSM-Thom polynomials} 
The SSM-Thom polynomial (Segre-Schwartz-\-Mac\-Pherson Thom polynomial) of a singularity, due to T. Ohmoto \cite{ohmotoSMTP}, plays the same role as its Thom polynomial, except the concept {\em fundamental class of closure} is replaced by {\em SSM class}. 

\smallskip

 Let $\eta\subset \E(m,m+\ell)$ be an orbit of the $\K(m,m+\ell)$ action. Consider the $G_{m,m+\ell}$ equivariant SSM class of $\eta$  
    \begin{multline}\label{eq:SSMeta}
    \ssm(\eta \subset \E(m,m+\ell)) \in H^{**}_{G_{m,m+\ell}}(\E(m,m+\ell))
    \\
    =H^{**}_{G_{m,m+\ell}}(pt)=\Q[[a_1,\ldots,a_m,b_1,\ldots,b_{m+\ell}]],
    \end{multline}
    where $a_i$ and $b_i$ are the same classes as in Section~\ref{sec:TPequivariantclass}. 

\begin{thm}[and Definition]\cite{ohmotoSMTP} 
\begin{itemize}
\item   The class in \eqref{eq:SSMeta} depends on $a_i$ and $b_i$ only through the quotient variables $c_i$ defined in~\eqref{eq:quotient_abstract}. When expressed as a formal power series in $c_i$ we call this formal power series the SSM-Thom polynomial of $\eta$, and denote it by $\T(\eta)$. 
\item We have
    \[ 
    T(\eta)=\Thom(\eta) + \text{higher degree terms}.
    \]
\item   If $\eta = \eta(Q,m,\ell)$ then $\T(\eta)$ only depends on $Q$ and $\ell$, not on $m$. We denote this power series by $\T(Q,\ell)$.
\item (Degeneracy locus interpretation.) For all maps $F : M^* \to N^{*+\ell}$ satisfying a
transversality condition we have
\begin{equation}
\label{eq:SSMTPdegeneracy}
\ssm( \eta(F) ) = T(\eta)(c(F)).
\end{equation}
\end{itemize}
\end{thm}

We hope the reader will pardon our use of the term `polynomial' to describe a power series (SSM-Thom). This choice is partly due to the term `Thom polynomial' becoming synonymous with the field and partly because, as we will demonstrate, our current calculations for these power series are limited to finite degrees. Additionally, `Thom series' refers to a different concept \cite{dstab}.

\begin{ex} \label{ex:A2l0}
The SSM-Thom polynomial $T(A_2,\ell=0)$ for the algebra $A_2=\C[x]/(x^3)$ and relative dimension $\ell=0$ is
    \begin{align*}
     &\ (c_{2}+ c_{1}^2) + (-3 c_{1}^3-3 c_{3}-6 c_{1} c_{2}) + (6 c_{1}^4+ 18 c_{1}^2 c_{2}+ 20 c_{1} c_{3}-c_{2}^2+ 7 c_{4}) \\
    & +(-76 c_{1}^2 c_{3}+ 6 c_{2} c_{3}+ 11 c_{1} c_{2}^2-40 c_{1}^3 c_{2}-56 c_{1} c_{4}-15 c_{5}-10 c_{1}^5)
    \\ 
    & + (-11 c_{2}^3+ 152 c_{1} c_{5}+ 31 c_{6}+ 15 c_{1}^6-53 c_{2} c_{4}+ 211 c_{1}^3 c_{3}+ 75 c_{1}^4 c_{2}
    \\ & \SK -39 c_{1} c_{2} c_{3}-46 c_{1}^2 c_{2}^2+ 22 c_{3}^2+ 245 c_{1}^2 c_{4}) \\
    & + (-400 c_{1} c_{6}+ 54 c_{1} c_{2}^3-21 c_{1}^7+ 199 c_{2} c_{5}-63 c_{7}-795 c_{1}^2 c_{5}-79 c_{3} c_{4}-126 c_{1}^5 c_{2} \\ & \SK + 423 c_{1} c_{2} c_{4}-199 c_{1} c_{3}^2-480 c_{1}^4 c_{3}+ 39 c_{2}^2 c_{3}-776 c_{1}^3 c_{4}+ 130 c_{1}^3 c_{2}^2+ 162 c_{1}^2 c_{2} c_{3}) \\ &
    + (28 c_{1}^8+ 39 c_{2}^4-119 c_{4}^2+ 127 c_{8}+ 953 c_{1}^5 c_{3}+ 907 c_{1}^2 c_{3}^2-17 c_{2} c_{3}^2+ 1997 c_{1}^4 c_{4}-93 c_{2}^2 c_{4} \\ & \SK + 1052 c_{1} c_{7}+ 196 c_{1}^6 c_{2}-295 c_{1}^4 c_{2}^2-155 c_{1}^2 c_{2}^3-512 c_{1}^3 c_{2} c_{3}-290 c_{1} c_{2}^2 c_{3} \\ & \SK -1855 c_{1}^2 c_{2} c_{4}+ 769 c_{1} c_{3} c_{4}-1755 c_{1} c_{2} c_{5}+ 2907 c_{1}^3 c_{5}+ 440 c_{3} c_{5}
    \\
    & \SK + 2468 c_{1}^2 c_{6}-742 c_{2} c_{6}) + \HOT.
\end{align*}
 We will describe the method that was used to calculate the displayed part of $T(A_2,\ell=0)$ in Section \ref{sec:IntThm}. However, it is important to note that the displayed degree 8 initial sum is the best we can do: we {do not know} how to calculate the higher degree parts of this particular power series. Several other SSM Thom polynomials (up to a certain degree) can be found on the Thom Polynomial Portal~\cite{TPP}.
\end{ex}


\subsection{The interpolation theorem} \label{sec:IntThm}
Let $\ell\geq 0$, the relative dimension, be fixed and let us choose a degree bound $d\leq M(\ell)$.

Let $Q_1, Q_2, \ldots, Q_r$ be the finite list of algebras that occur as local algebras of contact singularities of codimension $\leq d$, for germs $(\C^*,0)\to (\C^{*+\ell},0)$.  For $i=1,\ldots,r$ define
\begin{itemize}
\item $p_i:(\C^{s_i},0)\to (\C^{s_i+\ell},0)$ the prototype associated to $Q$ and $\ell$;
\item $G_i$ the maximal compact symmetry group (or its maximal torus) of $p_i$, with representations $\rhosou_i$ and $\rhotar_i$ in the source and the target;
\item $\psi_i:\Q[[c_1,c_2,\ldots]]\to H^{**}(BG_i)$ the homomorphism induced by
\[ 
1+c_1+c_2+\ldots \mapsto \frac{c(\rhotar_i)}{c(\rhosou_i)}.
\]
\end{itemize}
These concepts were discussed in Section~\ref{sec:geno_etc} with examples.
For a polynomial or power series $q$ we will write $q|_{\leq d}$ for the sum of its terms of degree up to and including $d$.

\begin{thm}[Interpolation theorem for SSM Thom polynomials] \label{thm:main} Let $d\leq M(\ell)$, and let $Q_i$ be the finite list of algebras that occur as local algebras of contact singularities of codimension $\leq d$ in $\E(*,*+\ell)$. The SSM Thom polynomial $T(Q,\ell)|_{\leq d}$ up to degree $d$ satisfies and is determined by the following axioms:
\begin{enumerate}
    \item $\psi_{i}( T(Q,\ell)|_{\leq d} ) = \left( e(\rhosou_i)/c(\rhosou_i) \right)|_{\leq d}$, if $Q=Q_i$;
    \item the degree $s_i,s_{i}+1,\ldots, d$ components of $\psi_i( T(Q,\ell)|_{\leq d} ) \cdot c(\rhosou_i)$ are 0, if $Q\not= Q_i$.
\end{enumerate}
\end{thm}

The proof will be presented in Sections~\ref{sec:proof1} and~\ref{sec:proof2}. Before that, we demonstrate how the theorem offers an efficient method for computing SSM-Thom polynomials.

\subsection{Example}
For $\ell=0$ let us choose $d=3\leq M(0)=8$. Then our list of algebras (see Figure~\ref{fig:l0}) is
\[
A_0, A_1, A_2, A_3.
\]
For each of these algebras (and $\ell=0$) we can carry out the definitions/algorithms of Section~\ref{sec:geno_etc} and we obtain: 

\noindent $\bullet$  $A_0$:
$e(\rhosou)=1, c(\rhosou)=1$, and $\psi$ is defined by
\[
\psi(1+c_1+c_2+\ldots)=\frac{1}{1}=1.
\]
That is, $\psi(c_1)=0$, $\psi(c_2)=0$, etc.

\noindent $\bullet$ $A_1$:
$e(\rhosou)=a, c(\rhosou)=1+a$, and $\psi$ is defined by
\[
\psi(1+c_1+c_2+\ldots)=\frac{1+2a}{1+a}=1+a-a^2+a^3-\ldots.
\]
That is, $\psi(c_1)=a$, $\psi(c_2)=-a^2$, etc.

\noindent $\bullet$ $A_2$:
$e(\rhosou)=2a^2, c(\rhosou)=(1+a)(1+2a)$, and $\psi$ is defined by
\[
\psi(1+c_1+c_2+\ldots)=\frac{1+3a}{1+a}=1+2a-2a^2+2a^3-\ldots.
\]

\noindent $\bullet$ $A_3$:
$e(\rhosou)=6a^3, c(\rhosou)=(1+a)(1+2a)(1+3a)$, and $\psi$ is defined by
\[
\psi(1+c_1+c_2+\ldots)=\frac{1+4a}{1+a}=1+3a-3a^2+3a^3-\ldots.
\]

Suppose we would like to calculate the SSM-Thom polynomial of $A_2$ up to degree $d=3$, that is, we want to determine the coefficients $v_\lambda$ in
\begin{multline}\label{eq:A1tpd3}
T(A_2,\ell=0)|_{\leq 3}=
v_0 + (v_1 c_1) + (v_{11}c_1^2+v_2 c_2)+
(v_{111}c_1^3+v_{12} c_1c_2 + v_3 c_3).
\end{multline}
According to Theorem~\ref{thm:main}, the constraints we obtain for our list of algebras for the coefficients $v_\lambda$ will uniquely determine them. Here are the constraints. 

\smallskip

\noindent $\bullet$ For $A_0$ we obtain that substituting $c_1=0, c_2=0,c_3=0$ in \eqref{eq:A1tpd3}, and then multiplying with $1$, the degree 0, 1, 2, 3 components are all 0. This imposes the constraint $v_0=0$. 

\noindent $\bullet$ For $A_1$ we obtain that substituting
$c_1=a, c_2=-a^2,c_3=a^3$ in \eqref{eq:A1tpd3}, and then multiplying with $(1+a)$, the degree 1, 2, 3 components are all 0. This imposes the constraints
\[
v_0+v_1=0, v_1-v_2+v_{11}=0, v_3-v_2+v_{11}-v_{12}+v_{111}=0.
\]

\noindent $\bullet$ For $A_3$ we obtain that substituting
$c_1=3a, c_2=-3a^2,c_3=3a^3$ in~\eqref{eq:A1tpd3}, and then multiplying with $(1+a)(1+2a)(1+3a)$, the degree 3 component is 0. This imposes the constraint
\[
6v_0+33v_1-18v_2+3v_3+54v_{11}-9v_{12}+27v_{111}=0.
\]

\noindent $\bullet$ Finally, for $A_2$ itself we obtain that substituting
$c_1=2a, c_2=-2a^2,c_3=2a^3$ in~\eqref{eq:A1tpd3}, up to degree 3 we must obtain 
\[
\frac{2a^2}{(1+a)(1+2a)}=2a^2-6a^3+\HOT.
\]
This imposes the constraints
\[
v_0=0, 2v_1=0, 4v_{11}-2v_2=2, 2v_3-4v_{12}+8v_{111}=-6.
\]

The constraints together have a unique solution (in accordance with the interpolation theorem), namely
\[
v_0=v_1=0, v_2=1, v_{11}=-6, v_3=-3, v_{12}=-6, v_{111}=-3,
\]
yielding
\[
T(A_2,\ell=0)|_{\leq 3}=(c_1^2+c_2)+(-3c_1^3-6c_1c_2-3c_3).
\]

\begin{rem}
    Of course, we know that the lowest degree part of an SSM Thom polynomial is the codimension of the singularity. For us $\codim(A_2)=2$, see Figure~\ref{fig:l0}, so we could have known that $v_0=v_1=0$ beforehand. We included these coefficients in the calculation to illustrate that the obvious geometric degree condition is a consequence of the interpolation theorem.
\end{rem}

If we chose the largest possible $d$ for $\ell=0$, namely $d=M(0)=8$, we would have obtained the degree 8 approximation of $T(A_2,\ell=0)$ given in Example~\ref{ex:A2l0}.

\begin{rem}
The degree 3 approximation we calculated above still falls within the range where the geometric method of \cite{ohmotoSMTP,nekarda1,nekarda2} remains applicable. What our calculations illustrate, even in this limited range, is that the computation is {\em automatic}: it does not require the understanding of the $A_2$ singularity loci near the prototype of all other singularities. That geometric task is generally too challenging; even deciding if that locus is empty or not is challenging (but cf. Section~\ref{sec:hierarchy}). In contrast, the degree 8 approximation goes far beyond the applicability of the geometric method.
\end{rem}

The next two sections contain the proof of Theorem~\ref{thm:main}.

\section{Part-I of the proof: Interpolation characterization\\ of CSM and SSM classes of orbits}
\label{sec:proof1}

In the first part of the proof we recall and improve the interpolation characterization of certain CSM classes \cite{RV,FRcsm} that is motivated by the Maulik-Okounkov axiom \cite{MO} for {\em stable envelopes}---a key concept in geometric representation theory.

\subsection{Interpolation characterization for CSM classes \cite{FRcsm}}
\label{sec:interpolation1}

Consider the linear representation $V$ of the algebraic group $G$. For each orbit $\omega\subset V$ let us choose a representative $x_\omega\in \omega$. We denote the stabilizer subgroup of $x_\omega$ by $G_\omega \subset G$, and let 
\[
\phi_{\omega} : H^*_G(V) \to H^*_G(\omega)
\]
be the restriction map. After the identifications $H^*_G(V)=H^*(BG)$, $H_G^*(\omega)=H^*_{G_\omega}(x_\omega)=H^*(BG_\omega)$ the map $\phi_\omega:H^*(BG)\to H^*(BG_\omega)$ is induced (through the functors $B$ and $H^*$) by the inclusion of groups $G_\omega \subset G$.

\begin{remark} \label{rem:polynomialring} \rm
    When $G$ is a linear algebraic group, such as the groups $G$ and $G_\omega$ above,
    then the rings $H^*(BG)=H^*(BG;\Q)$ are polynomial rings. This is proved in \cite{cartan} for compact connected groups; it holds for Lie groups because they are homotopy equivalent to their maximal compact subgroup; moreover, a finite quotient does not affect rational cohomology. In particular, these rings do not have zero divisors. We will use this fact throughout our proofs. 
\end{remark}

Let $T_\omega$ be the tangent space of $\omega$ at $x_\omega$ and let $N_\omega=V/T_\omega$ be the {normal} space. The group $G_\omega$ acts on $T_\omega$ and $N_\omega$, hence these representations have $G_\omega$ equivariant total Chern classes (denoted by $c(T_\omega), c(N_\omega)$), and equivariant Euler classes (denoted by $e(T_\omega), e(N_\omega)$). The product $c(T_\omega)c(N_\omega)$ is the $G_\omega$ equivariant total Chern class of $V$, or, in other words, $\phi_\omega(c(V))$.

\smallskip

\noindent We will make two assumptions on our representation.
\begin{enumerate}
    \item[(a)] The representation has finitely many orbits.
    \item[(b)] For every orbit $\omega$, the Euler class $e(N_\omega)\in H^*(BG_\omega)$ is not $0$ (equivalently, not a 0-divisor, cf. Remark~\ref{rem:polynomialring}).
\end{enumerate}

\begin{thm}[Thm 2.7 \cite{FRcsm}] \label{FRoriginal}
If a representation satisfies the conditions (a)-(b) above then the classes 
$\csm(\omega)\in H^*_G(V)$ satisfy and are determined by the properties (``axioms'')
\begin{enumerate}
    \item[(I'')]  $\phi_\omega(\csm(\omega))=c(T_\omega)e(N_\omega) \in H^*(BG_\omega)$;
    \item[(II'')]  $\phi_\theta(\csm(\omega))$ is divisible by $c(T_\omega)$ in $H^*(BG_\theta)$, for all orbits $\theta$;
    \item[(III'')]  $\deg(\phi_\theta(\csm(\omega)))< \deg(c(T_\theta)e(N_\theta))$ for all orbits $\theta\not=\omega$.
\end{enumerate}
\end{thm}

Recall that $\HH(X)\supset H^*(X)$ is the completion (with respect to cohomological degree). We will need a slightly improved version of the theorem above, namely:

\begin{thm} \label{FRimproved}
If a representation satisfies the conditions (a)-(b) above then the classes 
$\csm(\omega)\in \HH_G(V)$ satisfy and are determined by the properties (``axioms'')
\begin{enumerate}
    \item[(I')]  $\phi_\omega(\csm(\omega))=c(T_\omega)e(N_\omega) \in H^*(BG_\omega)$;
    \item[(II')]  $\phi_\theta(\csm(\omega))$ is in $H^*(BG_\theta)$, and in that ring it is divisible by $c(T_\omega)$, for all orbits $\theta$;
    \item[(III')]  $\deg(\phi_\theta(\csm(\omega)))< \deg(c(T_\theta)e(N_\theta))$ for all orbits $\theta\not=\omega$.
\end{enumerate}
\end{thm}

The proof of Theorem \ref{FRoriginal} given in \cite{FRcsm} serves as a proof of Theorem~\ref{FRimproved} with only cosmetic changes. 

\subsection{Interpolation characterization for SSM classes}
Now we are ready to state the {\em SSM version} of the theorems from Section \ref{sec:interpolation1}.

\begin{thm} \label{thm:SSM_interpolation}
    If a representation satisfies the conditions (a)-(b) above then the classes 
$\ssm(\omega)\in \HH_G(V)$ satisfy and are determined by the properties (``axioms'')
\begin{enumerate}
    \item[(I)]  $\phi_\omega(\ssm(\omega))=e(N_\omega)/c(N_\omega) \in \HH(BG_\omega)$;
    \item[(II)] if $\theta\not=\omega$, then $\phi_\theta(\ssm(\omega)) c(N_\theta)$
    \begin{enumerate}
       \item[(IIa)] is in $H^*(BG_\theta)$,
       \item[(IIb)] and is of degree $< \dim(N_\theta)$.
    \end{enumerate}
\end{enumerate}
\end{thm}

With the obvious extension of the notion {\em degree} $\in [0,\infty]$ from $H^*$ to $\HH$, of course, (IIb) implies (IIa). We included (IIa) for clarity.

\begin{proof}
   First we prove that $\ssm(\omega)$ satisfies the axioms. We have
\[
\phi_\omega(\ssm(\omega))
= \phi_\omega\left( \frac{\csm(\omega)}{c(V)} \right)
=  \frac{\phi_\omega(\csm(\omega))}{c(T_\omega)c(N_\omega)}
\overset{(I')}{=} \frac{c(T_\omega)e(N_\omega)}{c(T_\omega)c(N_\omega)}
= \frac{e(N_\omega)}{c(N_\omega)}
\]
proving (I). Now consider $\theta\not=\omega$ and
\begin{equation}\label{eq:temp}
\phi_\theta( \ssm(\omega)) c(N_\theta) 
= \frac{\phi_\theta(\csm(\omega))}{c(T_\theta)c(N_\theta)}c(N_\theta)
= \frac{\phi_\theta(\csm(\omega))}{c(T_\theta)}
\end{equation}
is in $H^*(BG_\theta)$ according to (II'), proving (IIa). Furthermore we obtain that 
\begin{multline*}
\deg\left( \phi_\theta( \ssm(\omega)) c(N_\theta)  \right)=
\deg\left( \phi_\theta(\csm(\omega)) \right) -
  \deg\left( c(T_\theta) \right) \\
\overset{(III')}{<} 
\deg(c(T_\theta)e(N_\theta)) - \deg\left( c(T_\theta) \right) 
\overset{(b)}{=} \deg\left( e(N_\theta) \right)
\overset{(b)}{=}\dim(N_\theta).
\end{multline*}
This proves (IIb). 

Now assume that a class $x(\omega)\in \HH_G(V)$ satisfies axioms (I)-(IIb), and define $y(\omega)=x(\omega)c(V)$. We will show that $y(\omega)$ satisfies (I')-(III'), which implies that $x(\omega)=\ssm(\omega)$.

Observe that 
\[
\phi_\omega(y(\omega))=
\phi_\omega(x(\omega) c(V)) 
\overset{(I)}{=}\frac{e(N_\omega)}{c(N_\omega)}\cdot c(N_\omega)c(T_\omega)
= e(N_\omega)c(T_\omega),
\]
proving (I'). We have
\[
\phi_\theta(y(\omega))
=\phi_\theta(x(\omega)c(V))
=\phi_\theta(x(\omega)) c(T_\theta) c(N_\theta)
\overset{(IIa)}{=} c(T_\theta) \times \text{polynomial},
\]
which proves (II'). 

We have
\begin{multline*}
\deg(\phi_\theta(y(\omega)))
=\deg(\phi_\theta(x(\omega)c(V))
\\
=\deg(\phi_\theta(x(\omega))c(T_\theta)c(N_\theta))
\leq \deg(\phi_\theta(x(\omega))c(N_\theta))+\deg(c(T_\theta))
\\
\overset{(IIb)}{<} \dim(N_\theta)+\deg(c(T_\theta))
\overset{(b)}{=} \deg( c(T_\theta)e(N_\theta)),
\end{multline*}
which proves (III'), and our proof is complete.
\end{proof}

\subsection{Interpolation characterization with degree bound}

Theorem~\ref{thm:SSM_interpolation} has a version in which the conditions and statements are modified to hold up to cohomological degree $2d$. These versions are useful in the theory of Thom polynomials, or other (finite dimensional approximations of) infinite dimensional representations. Hence in this section, we phrase and prove the needed ``up to degree $2d$'' version.

\medskip

Fix a positive integer $d$. Assume that 
\begin{enumerate}
    \item[($a^*$)] there is a $G$ invariant subvariety $\Pi\subset V$ of codimension~$>d$ such that  $V-\Pi$ has finitely many $G$-orbits, all of codimension $\leq d$;
    \item[($b^*$)] for every orbit $\omega\subset V-\Pi$, the Euler class $e(N_\omega)\in H^*(BG_\omega)$ is not 0. 
\end{enumerate}

We will consider cohomology classes up to (and including) degree $2d$. As before, we denote the truncation of a class $x$ up to (and including) cohomological degree $2d$ by $x|_{\leq d}$ or by $[x]_{\leq d}$. For an orbit $\omega\subset V-\Pi$ define $\ssmD(\omega)=\ssm(\omega)|_{\leq d}$.



\begin{thm}
Assume that the representation satisfies conditions ($a^*$)-($b^*$) above, and let $\omega\subset V-\Pi$ be an orbit. The class $\ssmD(\omega)\in \HD_G(V)$ satisfies and is determined by the properties (``axioms'')
\begin{enumerate}
    \item[(I${}^*$)]  $\phi_\omega(\ssmD(\omega))=\left[ e(N_\omega)/c(N_\omega)\right]_{\leq d} \in \HD(BG_\omega)$;
    \item[(II${}^*$)] $\deg\left( \left[ \phi_\theta(\ssmD(\omega)) c(N_\theta)\right]_{\leq d} \right)< \dim(N_\theta)$ for $\theta\not=\omega$, $\theta\subset V-\Pi$.
\end{enumerate}
\end{thm}

In (II${}^*$) we consider the class $\phi_\theta(\ssmD(\omega)) c(N_\theta)$. First we truncate it at degree $d$, then require that the components of degree $\geq \dim(N_\theta)$ are zero. This is equivalent to requiring that the degree $\dim(N_\theta), \dim(N_\theta)+1,\ldots,d$ components of $\phi_\theta(\ssmD(\omega)) c(N_\theta)$ vanish. 

\begin{proof}
We have 
\[
\phi_\theta(\ssmD(\omega))
=\phi_\theta([\ssm(\omega)]_{\leq d})
=[\phi_\theta(\ssm(\omega))]_{\leq d}
\overset{(I)}{=}[e(N_\omega)/c(N_\omega)]_{\leq d},
\]
which proves (I${}^*$). For $\theta\not=\omega$ we have
\[
\deg\left( 
[\phi_\theta(\ssmD(\omega)c(N_\theta)]_{\leq d}
\right)
=
\deg\left( 
[\phi_\theta(\ssm(\omega))c(N_\theta)]_{\leq d}
\right)
\overset{(II)}{<}
\dim(N_\theta),
\]
which proves (II${}^*$).

Now let us assume that two classes satisfy the properties (I${}^*$), (II${}^*$), and let $x(\omega)\in \HD_G(V)$ be their difference. For this class it holds that 
\begin{equation}\label{eq:temp_xx}
\deg\left( \left[ \phi_\theta(x(\omega)) c(N_\theta)\right]_{\leq d}\right)< \dim(N_\theta)
\end{equation}
for {\em every} $\theta\subset V-\Pi$ (both $\theta=\omega$ and $\theta\not=\omega$). Let $\theta_1,\theta_2,\ldots$ be a list of the orbits in $V-\Pi$ with $i<j \Rightarrow \theta_i\not\subset \overline{\theta}_j$. In particular, $\theta_1$ is the open orbit. 

We claim that $x(\omega)$ restricted to $\theta_1\cup \theta_2 \cup \ldots \cup \theta_s$ is 0 for all $s$. For $s=0$ the claim holds. Suppose we know this claim for $s-1$, that is, $x(\omega)$ restricted to $\theta_1\cup \theta_2 \cup \ldots \cup \theta_{s-1}$ is 0. It implies, via the usual Gysin sequence argument, see e.g. \cite[Lemma~2.6]{FRcsm}, that the class $\phi_{\theta_s}(x(\omega))$ is the $\leq d$ truncation of a class that is a multiple of $e(N_{\theta_s})$. The class $e(N_{\theta_s})$ is of degree $\dim(N_{\theta_s})$ (due to assumption~(b')), hence either 
$\phi_{\theta_s}(x(\omega))=0$ or its lowest degree non-0 term must have degree at least $\dim(N_{\theta_s})$. Multiplication by $c(N_\theta)$ ($=1+$higher degree terms) does not change the lowest degree term. Hence, according to  \eqref{eq:temp_xx} we must have $\phi_{\theta_s}(x(\omega))=0$. A Mayer-Vietoris argument, see e.g. Lemma~2.5 of \cite{FRcsm}, implies that $x(\omega)$ restricted to $\theta_1\cup \theta_2 \cup \ldots \cup \theta_s$ is zero, completing the induction step. 

We proved that $x(\omega)$ restricted to $V-\Pi$ is zero. However, the real codimension of $\Pi$ is at least $2d+2$, hence $\HD_G(V-\Pi)=\HD_G(V)$, and we obtain $x(\omega)=0\in \HD_G(V)$.
\end{proof}

\section{Part II of the proof}
\label{sec:proof2}

Choose $m\gg 0$ to satisfy $m\geq M(\ell)$ and that $(m+1)(m+\ell+1)>d$. Consider the representation of $G=\K(m,m+\ell)$ on the vector space $V=\E(m,m+\ell)$. As we reviewed in Section~\ref{sec:contact}, there exists a $G$-invariant subvariety $\Pi\subset V$ of codimension $>M(\ell)\geq d$ such that $V-\Pi$ has finitely many $G$-orbits, all of codimension $\leq d$. These are exactly the singularities $\eta(Q,m,\ell)$. That is, condition ($a^*$) from Section~\ref{thm:SSM_interpolation} is satisfied.

We claim that condition ($b^*$) also holds. We have no conceptual proof for this fact (cf. the discussion after Theorem~7.6 in \cite{MB}). In singularity theory this phenomenon is just considered to be a fact of life: {\em all singularities in $\E(m,m+\ell)-\Pi$ are quasi-homogeneous}, \cite[Thm.~7.6]{MB}. The way it is proved is that a $\C^*$-action is explicitly described for all of the finitely many ($<100$, but cf. Remark~\ref{rem:ListOfAlgebras}) algebras that may occur. Looking through the tables containing these algebras we can observe that in each case the $\C^*$-action is chosen in such a way that none of the weights on the source space are 0. Since the source space of a prototype is identified with the normal space of the singularity, it follows that the Euler class for this $\C^*$-action is a product of non-zero weights---that is, the Euler class is not 0. Therefore the Euler class for the whole symmetry group ($\supset \C^*$) cannot be zero either. This proves condition ($b^*$).

\begin{rem}
    We can be grateful to singularity theorists that they chose their $\C^*$ action that proves quasi-homogeneity in such a way that the Euler class in the source is not 0. They were not forced to do so. For example, continuing Example~\ref{ex:x2y3} we see that the Euler class of the action on the source space of the prototype for $(x^2,y^3)$, $\ell=0$ is 
    \[ ab(2a-b)(2b-a)(3b-a)2b(2b-a) \not= 0 \in \Q[a,b].\]
    If one was interested in quasi-homogeneity only---that is, in the existence of just one $\C^*$-action---then one could choose the first $\C^*$-component. This choice would result in $b=0$ and Euler class=0. Luckily, this is not the case. For all of the $Q$'s in the tables of \cite{MB} the choice of $\C^*$ is favorable.     
\end{rem}

Returning to the proof, let $\eta(Q_i,m,\ell)$ for $i=1,\ldots,r$ be the list of the finitely many orbits in $V-\Pi$. For each $i$ we use our previous notation: $G_i$ is the maximal compact symmetry group of the prototype $p_i:(\C^{s_i},0)\to (\C^{s_i+\ell},0)$, with its representations $\rhosou_i$ and $\rhotar_i$ in the domain and codomain. Recall that for each $i$ we have a homomorphism $\phi_i$ from 
\begin{multline*}
H^*_{\K(m,m+\ell)}(\E(m,m+\ell))
=H^*_{GL_m(\C)\times \GL_{m+\ell}(\C)}(pt) \\
=\Q[a_1,a_2,\ldots,a_m,b_1,b_2,\ldots,b_{m+\ell}]
\end{multline*}
to $H^*(BG_i)$. The map $\phi_i$ maps $a_i$ to the $i$'th Chern class of $\rhosou_i$, and $b_i$ to the $i$'th Chern class of $\rhotar_i$.

Recall that $Q$ is one of the $Q_i$, and hence $\eta(Q,m,\ell)$ is one of the orbits of $V-\Pi$. For the rest of the proof the following diagram will be useful.
\[
\begin{tikzcd}[column sep=10pt
    ,/tikz/column 2/.append style={anchor=base west}
    ]
T(Q,\ell)|_{\leq d}\ar[d,mapsto,"\kappa"] \in  & \Q[c_1,c_2,\ldots] \ar[d,"\kappa"] & & & H^*(BG_1) \\
\ssm_{\leq d}(\eta(Q,m,\ell)) \in & \Q[a_1,a_2,\ldots,a_m,b_1,b_2,\ldots,b_{m+\ell}] \ar[rrru,"\phi_1"]\ar[rrr,"\phi_2"]\ar[rrrdd,"\phi_r"]
& & & H^*(BG_2) \\
& & & & \vdots \\
& & & & H^*(BG_r). \\
\end{tikzcd}
\]
First we focus on the right side of the diagram. According to Theorem~\ref{thm:SSM_interpolation} the class 
$
S:=\ssmD(\eta(Q,m,\ell)) \in \Q[a_i,b_j]
$
is determined by the conditions
\begin{itemize}
\item $\phi_{i}(S)=\left[ e(\rhosou_i)/c(\rhosou_i) \right]_{\leq d}$ if $Q_i=Q$,
\item the degree $s_i,s_i+1,\ldots,d$ components of $\phi_{i}(S) c(\rhosou_i)$ vanish, if $Q_i\not=Q$.
\end{itemize}

Now we look at the top left part of the diagram. We claim that the map 
\[
\kappa=\kappa_{m,m+\ell}: \Q[c_1,c_2,\ldots] \to \Q[a_1,a_2,\ldots,a_m,b_1,b_2,\ldots,b_{m+\ell}]
\]
induced by 
\begin{equation}
 1+c_1+c_2+\ldots=\frac{1+b_1+b_2+\ldots+b_{m+\ell}}{1+a_1+a_2+\ldots+a_m}
\end{equation}
is injective up to degree $d$. 

Indeed, this map is studied in algebraic combinatorics where it is proved that 
\begin{equation}\label{eq:ker}
\ker(\kappa_{m,m+\ell})= \text{span}\{ s_\lambda : \lambda_{m+1} \geq m+\ell+1\},
\end{equation}
where 
\begin{equation}\label{eq:Schur}
s_\lambda=s_{(\lambda_1\geq \lambda_2\geq \ldots \geq \lambda_n)}=
\det\begin{pmatrix}
    c_{\lambda_i+j-i}
\end{pmatrix}_{i,j=1,\ldots,n}
\end{equation}
is the Schur polynomial associated with the partition $\lambda$. 

The statement~\eqref{eq:ker} is proved (albeit not stated) in \cite[Section 3.2]{FP}.  Indeed, the $\supset$ part follows from the factorization formula on page 37 of that book, and the $\subset$ part follows from the Proposition on page 36 of \cite{FP}---cf. \cite[Section 6.1]{FRcsm}.

The $\lambda_{m+1}\geq m+\ell+1$ condition forces that all the generators of the kernel is of degree at least $(m+1)(m+\ell+1)>d$. Therefore, indeed, $\kappa_{m,m+\ell}$ is injective up to (and including) degree $d$.

Putting together our two arguments we obtained that $T=T(Q,\ell)|_{\leq d}\in \Q[c_i]$ is determined by the conditions
\begin{itemize}
\item $\phi_{i}\kappa(T)=\left[ e(\rhosou_i)/c(\rhosou_i) \right]_{\leq d}$ if $Q_i=Q$,
\item the degree $s_i,s_i+1,\ldots,d$ components of $\phi_{i}\kappa(T) c(\rhosou_i)$ vanish, if $Q_i\not=Q$.
\end{itemize}
Since $\phi_i\kappa=\psi_i$, the proof is complete. \qed

\section{Expansions, algebraic structure}\label{sec:expansions}

\subsection{Expansions}
Thom polynomials and various other polynomials in enumerative geometry show interesting properties when expanded in appropriate bases. These bases include the Chern monomial basis, the Schur basis, and the Schur-tilde basis.

The Schur functions are defined in \eqref{eq:Schur}. Schur tilde functions (in fact, formal power series) are defined in \cite[Def.~8.2]{FRcsm}:
\[
\tilde{s}_\lambda=
S \left(
\prod_{i=1}^k \left( \frac{z_i}{1+z_i}\right)^{\lambda_i}
\prod_{j=1}^{\infty}
\prod_{i=1}^j
\frac{1+z_i-z_j}{1+z_i}
\right)
\]
where $\lambda$ is a partition, and the the linear $S$ operation turns a $z$-monomial 
$z_1^{\mu_1}\cdots z_k^{\mu_k}$ to 
\[
\det(c_{\mu_i+j-i})_{i,j=1,\ldots,k}.
\]
Notice that this last determinant is seemingly the same formula as the definition of the Schur function. However, here the integer vector $\mu$ is not necessarily a partition.
For example 
\[
\tilde{s}_{41}=
s_{41}-(3s_{411}+3s_{42}+5s_{51})+
(6s_{4111}+10s_{42}+5s_{43}+16s_{511}+16s_{52}+15s_{61})-\HOT.
\]

While Schur functions are equivariant fundamental classes of matrix Schubert varieties, Schur-tilde functions are (suitable limits of) equivariant SSM classes of the same. Hence, it is natural to expect that formulas for SSM classes of geometrically relevant varieties should have favorable properties in their Schur-tilde expansions. Evidence for this expectation is provided in \cite{FRcsm, sutipoj}.

\begin{remark} \rm
    One reality check for SSM formulas stems from the additivity property of SSM classes (cf. (iv) of Section~\ref{sec:CSM}). It implies that the sum of the SSM-Thom polynomials of {\em all} singularities is the SSM class of the whole ambient space, which is by definition $=1$. In Chern monomial basis and in Schur basis this means that in the sum of all SSM-Thom polynomials all $c_\lambda$, $s_\lambda$ terms must vanish, except $c_0=s_0=1$. In Schur-tilde basis this means that the sum of all SSM-Thom polynomials is $\sum_\lambda \tilde{s}_\lambda$. (Schur-tilde functions act like a formal power series-valued probability distribution on partitions.) Verifying the phenomena of this remark in Figures~\ref{fig:S0} and~\ref{fig:S1} is instructive. 
\end{remark}

\subsection{Simplification}
Notable simplifications, such as fewer terms and smaller absolute value coefficients, can be observed in the Schur and Schur-tilde expansions of SSM-Thom polynomials. The reader is encouraged to examine these patterns in Figures~\ref{fig:S0} and~\ref{fig:S1}, as well as in the extensive examples found on the \cite{TPP}.
 
\begin{figure}
\[
\begin{adjustbox}{angle=90}
\begin{tabular}{| C | C | C | C | C | C | C |}
\hline
& \deg 0 & \deg 1 & \deg 2 & \deg 3 & \deg 4 & \deg 5 \\
\hline
{A_0}&1&-s_{1}&s_{2}+s_{11}&-s_{3}-2\,s_{21}-s_{111}&s_{4}+s_{22}+3\,s_{31}+3\,s_{211}+s_{1111}&-s_{5}-2\,s_
{32}-4\,s_{41}-2\,s_{221}\hfill\ 
\\ 
 & & & & & & \ \hfill -6\,s_{311}-4\,s_{2111}-s_{11111}
\\ \hline
{A_1}&&s_{1}&-3\,s_{2}-2\,s_{11}&7\,s_{3}+9\,s_{21}+3\,s_{111}&-15\,s_{4}-11\,s_{22}-28\,s_{31} \hfill\  
&31\,s_{5}+48\,s_{32}+75\,s_{41}+32\,s_{221}\hfill\ 
\\ 
 & & & & & \ \hfill -18\,s_{211}-4\,s_{1111} & \ \hfill +70\,s_{311}+30\,s_{2111}+5\,s_{11111}
\\
\hline
{A_2}&&&2\,s_{2}+s_{11}&-12\,s_{3
}-12\,s_{21}-3\,s_{111}&50\,s_{4}+29\,s_{22}+73\,s_{31} \hfill\ &-180\,s_{5}-218\,s_{32}-340\,s_{41}-119\,s_{221}\hfill\ 
\\ 
 & & & & & \ \hfill +36\,s_{211}+6\,s_{1111} & \ \hfill -245\,s_{311}-80\,s_{2111}-10\,s_{11111}
\\
\hline
{A_3}&&&&6\,s_{3}+5\,s_{21}+s_{111}&-60\,s_{4}-30\,s_{22}-74\,s_{31} \hfill\ 
&390\,s_{5}+398\,s_{32}+625\,s_{41}+189\,s_{221}\hfill\ 
\\ 
 & & & & & \ \hfill -30\,s_{211}-4\,s_{1111}& \ \hfill +375\,s_{311}+100\,s_{2111}+10\,s_{11111}
\\
\hline
{A_4}&&&&&24\,s_{4}+10\,s_{22}+26\,s_{31} \hfill\ &-360\,s_{5}-316\,s_{32}-510\,s_{41}-133\,s_{221}\hfill\ 
\\
& & & & & \ \hfill +9\,s_{211}+s_{1111} & \ \hfill -265\,s_{311}-60\,s_{2111}-5\,s_{11111}
\\
\hline
{A_5}&&&&&&120\,s_{5}+92\,s_{32}+154\,s_{41}+35\,s_{221}\hfill\ 
\\ 
 & & & & & & \ \hfill +71\,s_{311}+14\,s_{2111
}+s_{11111}
\\
\hline
{I_{22}}&&&&&s_{22}&-6\,s_{32}-4\,s_{221}
\\
&&&&&&
\\
\hline
{I_{23}}&&&&&&4\,s_{32}+2\,s_{221}
\\
&&&&&&
\\
\hline
 \\
 \\
\hline
{A_0}&\tilde{s}_{0} &&&&&
\\
&&&&&&
\\ \hline
{A_1}&&\tilde{s}_{1} &-
\tilde{s}_{2} &\tilde{s}_{3} +2\,\tilde{s}
_{21} &-\tilde{s}_{4} -4\,\tilde{s}
_{31} &\tilde{s}_{5} +6\,\tilde{s}
_{41} +4\,\tilde{s}_{311} 
\\
&&&&&&
\\ \hline
{A_2}&&&2\,\tilde{s}_{2} +\tilde{s}_{11} &-6\,\tilde{s}_{3} -6\,\tilde{s}_{21} &14\,\tilde{s}_{4} +26\,
\tilde{s}_{31} +5\,\tilde{s}_{22} +5\,\tilde{s}_{211} &-30\,\tilde{s}_{5} -82\,
\tilde{s}_{41} -34\,\tilde{s}_{32}\hfill\ 
\\ 
 & & & & & & \ \hfill  -38\,
\tilde{s}_{311} -12\,\tilde{s}_{221} 
\\ \hline
{A_3}&&&&6\,\tilde{s}_{3} +5
\,\tilde{s}_{21} +\tilde{s}_{111} &-36
\,\tilde{s}_{4} -47\,\tilde{s}_{31} -15\,
\tilde{s}_{22} -13\,\tilde{s}_{211} &
150\,\tilde{s}_{5} +281\,\tilde{s}_{41} +
150\,\tilde{s}_{32}\hfill\ 
\\ 
 & & & & & & \ \hfill  +124\,\tilde{s}_{311} +53\,\tilde{s}_{221} +9\,\tilde{s}_{2111} 
\\ \hline
{A_4}&&&&&24\,\tilde{s}
_{4} +26\,\tilde{s}_{31} +10\,\tilde{s}
_{22} +9\,\tilde{s}_{211} +\tilde{s}
_{1111} &-240\,\tilde{s}_{5} -358\,\tilde{s}_{41} -208\,\tilde{s}_{32}\hfill\ 
\\ 
 & & & & & & \ \hfill  -160
\,\tilde{s}_{311} -76\,\tilde{s}_{221} -22\,\tilde{s}_{2111} 
\\ \hline
 {A_5}&&&&&&120\,\tilde{s}_{5} +154\,\tilde{s}_{41} +92\,\tilde{s}_{32}\hfill\ 
\\ 
 & & & & & & \ \hfill  +71\,\tilde{s}_{311} +35
\,\tilde{s}_{221} +14\,\tilde{s}_{2111} +\tilde{s}_{11111} 
 \\ \hline
{I_{22}}&&&&&\tilde{s}_{22} &-3\,\tilde{s}_{32} -\tilde{s}_{221} 
\\
&&&&&&
\\ \hline
{I_{23}}&&&&&&4\,\tilde{s}_{32} +2\,\tilde{s}
_{221} 
\\
&&&&&&
\\
\hline
\end{tabular} 
\end{adjustbox}
\]
\caption{SSM-Thom polynomials for $l=0$ up to degree 5}\label{fig:S0}
\end{figure}

\begin{figure}
\[
\begin{adjustbox}{angle=90}
\begin{tabular}{| C | C | C | C | C | C | C | C|}
\hline
& \deg 0 & \deg 1 & \deg 2 & \deg 3 & \deg 4 & \deg 5 & \deg 6\\
\hline
{A_0}&s_{0}&0&-s_{2}&2\,s_{3}+s_{21}&-3\,s_{4}-3\,s_{31}-s_{211}
&4\,s_{5}+6\,s_{41}
&-5\,s_{6}+s_{33}-10\,s_{51} 
\\ & & & & & & +4\,s_{311}+s_{2111}  &
-10\,s_{411}-5\,s_{3111}-s_{21111}
\\ \hline
{A_1}&&&s_{2}&-2\,s_{3}-s_{21}&-s_{4}-s_{22}+s_{31}+s_{211}&20\,s_{5}+10\,s_{32}+14\,s_{41}
&-
87\,s_{6}-20\,s_{33}-57\,s_{42}
\\ & & & & & & +2\,s_{221}-s_{2111} &
-96\,s_{51}-s_{222}-23\,s_{321}
\\ & & & & & & &
-32\,s_{411}-3\,s_{2211}-s_{3111}+s_{21111}
\\ \hline
{A_2}&&&&&4\,s_{
4}+s_{22}+2\,s_{31}&-24\,s_{5}-10\,s_{32}-20\,s_{41}&56\,s_{6}+13\,s_{33}+38\,s_{42}
\\ & & & & & & -2\,s_{221}-4\,s_{311} &
+76\,s_{51}+18\,s_{321}+36\,s_{411}
\\ & & & & & & &
+3\,s_{2211}+6\,s_{3111}
\\ \hline
{A_3}&&
&&&&
&36\,s_{6}+5\,s_{33}+19\,s_{42}+30\,s_{51}+
\\ & & & & & & &
s_{222}+5\,s_{321}+6\,s_{411}
\\ 
\hline
{III_{22}}&&&&&&&s_{33}
\\
&&&&&&&
\\
\hline
\\
\\
\\
\hline
{A_0}&\tilde{s}_{0} &
\tilde{s}_{1} &\tilde{s}_{11} &\tilde{s}_{111} &\tilde{s}_{1111} &\tilde{s}_{11111} &\tilde{s}_{111111} 
\\ \hline
{A_1}&&&\tilde{s}_{2} &
\tilde{s}_{3} +\tilde{s}_{21} &-3\,\tilde{s}_{4} -\tilde{s}_{31} +\tilde{s}_{211} &5\,\tilde{s}_{5} +5\,\tilde{s}_{41}  
&-7\,\tilde{s}_{6} -13\,\tilde{s}_{51} +\tilde{s}_{411} 
\\ & & & & & &+2\,\tilde{s}_{32} -\tilde{s}_{311} +\tilde{s}_{2111} &
+2\,\tilde{s}_{321} -\tilde{s}_{3111} +\tilde{s}_{21111} 
\\ \hline
{A_2}&&&&&4\,\tilde{s}_{4} +2\,\tilde{s}_{31} +\tilde{s}_{22} &-4\,\tilde{s}_{5} -4\,\tilde{s}_{41} 
&-28\,\tilde{s}_{6} -16\,\tilde{s}_{51} -18\,\tilde{s}_{42}
\\ & & & & & & -\tilde{s}_{32} +2\,\tilde{s}_{311} +\tilde{s}_{221}  &
-5\,\tilde{s}_{33} -6\,\tilde{s}_{411} -6\,\tilde{s}_{321}
\\ & & & & & & &
+2\,\tilde{s}_{3111} +\tilde{s}_{2211} 
\\ \hline
{A_3}&&&&&&
&36\,\tilde{s}_{6} +30\,\tilde{s}_{51} +19\,\tilde{s}_{42}
+5\,\tilde{s}_{33} 
\\ & & & & & & &
+6\,\tilde{s}_{411} 
+5\,\tilde{s}_{321} +\tilde{s}_{222} 
\\ \hline
{III_{22}}&&&&&&&\tilde{s}_{33} 
\\
&&&&&&&
\\
\hline
\end{tabular} 
\end{adjustbox}
\]
\caption{SSM-Thom polynomials for $l=1$ up to degree 6}\label{fig:S1}
\end{figure}

One spectacular instance of simplification is that while the Chern monomial and Schur expansions of $T(A_0,\ell)$ are complicated looking expressions, we have
\[
T(A_0,\ell=0)=\tilde{s}_0, \qquad \text{a one-term expression,}
\]
and more generally
\[
T(A_0,\ell)=\sum_{\lambda: \lambda_1\leq\ell} \tilde{s}_\lambda.
\]

\subsection{Positivity, stabilization}

Ordinary (not SSM-) Thom polynomials 
\begin{itemize}
\item are Schur positive \cite{pragacz:positivity},
\item are conjectured to be monomial positive if and only if the singularity is $A_n$ \cite{rrtp}, and
\item satisfy the so-called $d$-stability property, related to increasing the relative dimension $\ell$ \cite{dstab}.
\end{itemize}

A large number of (initial sums of) SSM-Thom polynomials that we were able to calculate support the following conjecture.

\begin{con}
For $\ell=0$ the SSM-Thom polynomials have alternating signs both in Schur and Schur-tilde expansions.
\end{con}

Other, more mysterious, sign rules can also be observed. Apparently, the Schur expansions of $T(Q,\ell=1)$ have positive coefficients for many $Q$, for example $Q=A_4,A_5,I_{22},\ldots$, but not for $Q=A_1,A_2,A_3$.

We do not know how the mentioned d-stability property of Thom polynomials generalizes to SSM-Thom polynomials.

\section{Applications} \label{sec:app}

We believe that SSM-Thom polynomials will be effective tools in enumerative geometry, differential topology, and singularity theory itself. In this section, we outline three promising avenues through which SSM-Thom polynomials could find applications. We provide concrete examples to illustrate each directions and leave the systematic investigation of such applications for future research.

\subsection{Enumerative geometry}\label{sec:application}

Consider a degree $d$ map $F:\PP^5 \to \PP^6$ and assume it is stable. A glance at Figure~\ref{fig:l1} tells us that it has only $A_0$, $A_1$ and $A_2$ singularities, and that the $A_2$ singularity locus $A_2(F)$ is a smooth complex curve.

\begin{thm}
    We have
    \begin{align*}
    \deg( A_2(F) \subset \PP^5)  & =  21(d-1)^2(6d-7)^2 \\
    \chi(A_2(F) ) & = -14(d-1)^2(1048d^3-3928d^2+4887d-2019).
    \end{align*}
\end{thm}

The first statement is what the ordinary (not SSM-) Thom polynomial would have given us. The second statement offers a glimpse into the additional geometric information encoded by SSM-Thom polynomials.

\begin{proof}
Let $h\in H^*(\PP^5)$ be the hyperplane class. Then 
\begin{multline*}
    c(F)=\frac{F^*(c(T\PP^6))}{c(T\PP^5)}=
    \frac{(1+dh)^7}{(1+h)^6}=\\
\underbrace{(7d-6)h}_{c_1(F)} + \underbrace{21(d-1)^2 h^2}_{c_2(F)}+ \underbrace{7(d-1)^2(5d-8)h^3}_{c_3(F)}+ \\
\underbrace{7(d - 1)^2(5d^2 - 20d + 18)h^4}_{c_4(F)}
+\underbrace{21(d - 1)^2(d - 2)(d^2 - 6d + 6)h^5}_{c_5(F)}.
\end{multline*}
Substituting these classes in $T(A_2,\ell=1)|_{\leq 5}=$
\[
(c_1c_3 + c_2^2 + 2c_4) + (-2c_1^2c_3 - 2c_1c_2^2 - 8c_1c_4 - 4c_2c_3 - 8c_5),
\]
we obtain that 
\begin{multline*}
    \ssm(A_2(F)\subset \PP^5)=\\
21(d - 1)^2(6d - 7)^2h^4 - 14(d - 1)^2(1048d^3 - 3604d^2 + 4131d - 1578) h^5.
\end{multline*}
Using $\ssm(A_2(F)\subset \PP^5)=\csm(A_2(F)\subset \PP^5)/c(T\PP^5)$ we have
\begin{multline*}
    \csm(A_2(F)\subset \PP^5)=\\
21(d - 1)^2(6d - 7)^2 h^4 - 14(d - 1)^2(1048d^3 - 3928d^2 + 4887d - 2019)h^5.
\end{multline*}
The Ohmoto-Aluffi theorem describes the algebra of how the coefficients of CSM classes encode the Euler characteristics of (general) linear sections of a projective variety.

\begin{thm}\cite{ohmotoCHI, aluffi} 
    For $\csm(\Sigma \subset \PP^N)=\sum a_i h^i$ define $\gamma_{\Sigma}(t)=\sum_i a_i t^{N-i}$ and 
    \[
    \chi_{\Sigma}(t)=\sum_r \chi(\ \  \Sigma \cap\underbrace{ H_1 \cap \ldots \cap H_r}_{\text{$r$ general hyperplanes}})(-t)^r,
    \]
    Then the polynomials $\gamma_{\Sigma}(t)$ and $\chi_{\Sigma}(t)$ are related by the involution $p \leftrightarrow \mathcal{A}(p)$,
    \[
    \mathcal{A}(p)(t)=\frac{tp(-t-1)+p(0)}{t+1}.
    \]
\end{thm}

Applying this theorem to our situation gives
\begin{multline*}
    \chi_{A_2(F)}(t)=\\
-14(d - 1)^2(1048d^3 - 3928d^2 + 4887d - 2019) - 21(d - 1)^2(6d - 7)^2t,
\end{multline*}
and the proof is complete.
\end{proof}

\subsection{Hierarchy of singularities}\label{sec:hierarchy}
In this section we illustrate the results of \cite[\S 6]{FeherPatakfalvi} by analyzing a specific example involving the hierarchy of singularities for $\ell=1$, as depicted in Figure~\ref{fig:l1}. For instance, the figure shows that the singularity defined by the algebra 
\[
Q=\C[x,y,z]/(x^2+y^2+z^2,xy,yz,zx)
\]
lies in the closure of the singularity defined by $I_{24}=\C[x,y]/(xy,x^2+y^4)$, but not in the closure of the singularity defined by $I_{33}=\C[x,y]/(xy,x^3+y^3)$. This distinction underscores the nuanced structure of adjacency relations within the singularity hierarchy (note that both $I_{24}$ and $I_{33}$ are of Thom-Boardman type $\Sigma^{20}$, that is, they cannot be distinguished by arbitrary high derivatives.)

To prove this statement directly from the definition, one must analyze the local algebras of all germs near $0$ in the prototype of $Q$ for $\ell=1$. This involves finding local algebras arbitrarily close to $0$ that are isomorphic to $I_{24}$, while also showing the existence of a neighborhood of $0$ that contains no algebras isomorphic to $I_{33}$. Determining whether a given algebra (defined by presentation) is isomorphic to $I_{24}$ or $I_{33}$ is a challenging problem in commutative algebra, requiring a detailed analysis of invariants and structural properties.

Thom polynomials can help us decide which algebras are comparable in the hierarchy without commutative algebra. Consider the homomorphism 
\[
\psi_Q:\Q[c_1,c_2,\ldots]\to \Q[a]
\] 
associated with the prototype of $Q$, $\ell=1$, namely the one induced by 
\[
1+c_1+c_2+\ldots \mapsto \frac{ (1+2a)^4}{(1+a)^3},
\]
that is, $c\mapsto 5a$, $c_2\mapsto 6a^2$, $c_3\mapsto -2a^3$, etc. The SSM-Thom polynomials up to degree 14, in particular the Thom polynomials $\Thom(I_{24},\ell=1)$, $\Thom(I_{33},\ell=1)$ can be computed with the interpolation theorem. We have $\Thom(I_{24},\ell=1)=$
\begin{multline*}
-6 c_1^3 c_2 c_6+ 9 c_1^3 c_3 c_5-3 c_1^3 c_4^2-5 c_1^2 c_2^2 c_5+ 2 c_1^2 c_2 c_3 c_4+ 3 c_1^2 c_3^3-2 c_1 c_2^3 c_4+ 2 c_1 c_2^2 c_3^2\\
-54 c_1^2 c_2 c_7+ 67 c_1^2 c_3 c_6-13 c_1^2 c_4 c_5-33 c_1 c_2^2 c_6+ 19 c_1 c_2 c_3 c_5-7 c_1 c_2 c_4^2+ 21 c_1 c_3^2 c_4\\
-6 c_2^3 c_5+ 4 c_2^2 c_3 c_4+ 2 c_2 c_3^3-156 c_1 c_2 c_8+ 162 c_1 c_3 c_7+ 24 c_1 c_4 c_6-30 c_1 c_5^2-52 c_2^2 c_7\\+ 28 c_2 c_3 c_6-12 c_2 c_4 c_5+ 16 c_3^2 c_5+ 20 c_3 c_4^2-144 c_2 c_9+ 128 c_3 c_8+ 52 c_4 c_7-36 c_5 c_6,
\end{multline*}
and the reader can find the other one on the Thom polynomial portal \cite{TPP}. Substitutions give:
\[
\psi_Q (\Thom(I_{24},\ell=1)) = 6a^{11}, \qquad \psi_Q(\Thom(I_{33},\ell=1))=0.
\]
The geometric meaning of $\psi_Q (\Thom(I_{24}))$ is the $\C^*$ equivariant fundamental class of the $I_{22}$ singularity locus in the domain of the prototype for $Q,\ell=1$. Hence, the very fact that it is not 0 implies that that locus is not empty: $Q$ is below $I_{22}$ in the hierarchy, see Figure~\ref{fig:l1}.

The fact that $\psi_Q(\Thom(I_{33},\ell=1))=0$ does not {\em a priori} imply that the $I_{33}$ locus near 0 in the domain of the prototype of $Q,\ell=1$ is empty. However, it is proved in \cite{FeherPatakfalvi} that under a {\em positivity} condition on the weights of $\rhosou$ for $Q$, it does. This positivity condition obviously holds in our case, proving that---as depicted in Figure~\ref{fig:l1}---$Q$ is not below $I_{33}$ in the hierarchy of singularities. 

\subsection{Complex vs real forms of singularities} \label{sec:RC}

Our method of calculating characteristic classes of singularities, namely Thom polynomials and SSM-Thom polynomials, also has consequences about the classification of complex singularities in Mather's nice dimensions. We are going to illustrate this with an example.

In \cite{MB} a list of {\em real} contact singularities is given, in particular the two algebras listed in lines 1 and 2 in Table~7.8 are:
\[
Q_1=
\R[x,y,z]/(x^2+y^2+z^2, xy,xz,yz),
\quad
Q_2=
\R[x,y,z]/(x^2,y^2,z^2-xy,xz+yz).
\]
We claim that the complex versions (replace $\R$ with $\C$) of $Q_1$ and $Q_2$ are isomorphic. Theoretically this statement could be proved by finding the concrete isomorphism, our point is that it follows from Thom polynomials and their interpolation properties.

Namely, consider the obvious $U(1)$ symmetry of both prototypes (say for $\ell=1$):
\[
\rhosou=12\lambda, \qquad  \rhotar=4\lambda^2 \oplus 9\lambda.
\]
Since they are the same, the homomorphisms $\psi_{Q_1},\psi_{Q_2}:\Q[c_1,c_2,
\ldots]\to \Q[a]$ associated with these two algebras are the same, they are both induced by 
\[
1+c_1+c_2+\ldots \mapsto \frac{(1+2a)^4}{(1+a)^3}.
\]
Moreover, part $(I^*)$ of the interpolation theorem gives 
\[
\psi_{Q_1} ( \Thom(Q_1,\ell) ) =a^{12}.
\]
If the two algebras were not isomorphic, then part $(II^*)$ of the interpolation theorem would give  
\[
\psi_{Q_2} (\Thom(Q_1,\ell) )=0.
\]
Since $\psi_1=\psi_2$ these two conditions contradict. The complex forms of the algebras $Q_1$ and $Q_2$---or any pair of algebras with the same symmetry and representations---are isomorphic.

\bibliographystyle{alpha}
\bibliography{references}

\begin{thebibliography}{AMSS23}

\bibitem[Alu13]{aluffi}
P.~Aluffi.
\newblock Euler characteristics of general linear sections and polynomial {C}hern classes.
\newblock {\em Rend. Circ. Mat. Palermo (special issue)}, pages 3--26, 2013.

\bibitem[AMSS23]{AMSS_csm}
P.~Aluffi, L.~Mihalcea, J.~Sch\"urmann, and Ch. Su.
\newblock Shadows of characteristic cycles, {V}erma modules, and positivity of {C}hern-{S}chwartz-{M}ac{P}herson classes of {S}chubert cells.
\newblock {\em Duke Math. J.}, 172(17):3257--3320, 2023.

\bibitem[B{\'e}r20]{BercziNew}
G.~B{\'e}rczi.
\newblock Non-reductive geometric invariant theory and {T}hom polynomials.
\newblock arXiv:\-2012.06425, 2020.

\bibitem[BS81]{BrSch}
J.-P. Brasselet and M.-H. Schwartz.
\newblock Sur les classes de chern d’un ensemble analytique complexe.
\newblock {\em Ast\`erisque}, 83:93--147, 1981.
\newblock Soc. Math. France.

\bibitem[BS12]{BSz}
G.~B\'erczi and A.~{Sz}enes.
\newblock Thom polynomials of {M}orin singularities.
\newblock {\em Annals of Math.}, 175:567--–629, 2012.

\bibitem[Car51]{cartan}
H.~Cartan.
\newblock La transgression dans un groupe de lie et dans un espace fibr\'e principal.
\newblock In {\em Colloque de topologie (Bruxelles, 1950), Georges Thone, Li\`ege \& Masson, Paris}, 1951.

\bibitem[Dam72]{damon}
J.~Damon.
\newblock Thom polynomials for contact class singularities.
\newblock Ph.D. Thesis, Harvard, 1972.

\bibitem[FP98]{FP}
W.~Fulton and P.~Pragacz.
\newblock {\em Schubert varieties and degeneracy loci}.
\newblock Number 1689 in LNM. Springer, 1998.

\bibitem[FP09]{FeherPatakfalvi}
L.~Feh\'er and Zs. Patakfalvi.
\newblock The incidence class and the hierarchy of orbits.
\newblock {\em Open Mathematics}, 7(3):429--441, 2009.

\bibitem[FR07]{dstab}
L.~M. Feh{\'e}r and R.~Rim{\'a}nyi.
\newblock On the structure of {T}hom polynomials of singularities.
\newblock {\em Bull. London Math. Soc.}, 39:541--549, 2007.

\bibitem[FR12]{FRannals}
L.~M. Feh{\'e}r and R.~Rim{\'a}nyi.
\newblock Thom series of contact singularities.
\newblock {\em Annals of Math.}, 176(3):1381--1426, 2012.

\bibitem[FR18]{FRcsm}
L.~M. Feh\'er and R.~Rim\'anyi.
\newblock Chern-{S}chwartz-{M}ac{P}herson classes of degeneracy loci.
\newblock {\em Geometry and Topology}, 22:3575--3622, 2018.

\bibitem[Kaz]{kaza:noas}
M.~\'E. Kazarian.
\newblock Non-associative {H}ilbert scheme and {T}hom polynomials.
\newblock unpublished.

\bibitem[Mac74]{macpherson}
R.~MacPherson.
\newblock Chern classes for singular algebraic varieties.
\newblock {\em Ann. of Math.}, 100:421--432, 1974.

\bibitem[Mat69]{mather4}
John~N. Mather.
\newblock Stability of {$C\sp{\infty }$} mappings. {IV. C}lassification of stable germs by {$R$}-algebras.
\newblock {\em Inst. Hautes \'Etudes Sci. Publ. Math.}, 37:223--248, 1969.

\bibitem[Mat71]{mather6}
J.~N. Mather.
\newblock Stability of {${C}^{\infty}$} mappings {VI}: the nice dimensions.
\newblock In C.~T.~C. Wall, editor, {\em Proceedings of Liverpool Singularities I}, LNM 192. Springer, 1971.

\bibitem[MNB20]{MB}
D.~Mond and J.~J. {Nu\~no}-Ballesteros.
\newblock {\em Singularities of mappings}.
\newblock Number 357 in Grundlehren der mathematischen Wissenschaften. Springer, 2020.

\bibitem[MO19]{MO}
D.~Maulik and A.~Okounkov.
\newblock {\em Quantum Groups and Quantum Cohomology}, volume 408 of {\em Ast\'erisque}.
\newblock SMF, 2019.

\bibitem[NOa]{nekarda1}
S.~Nekarda and T.~Ohmoto.
\newblock Computing higher {T}hom polynomials for multi-singularity classes.
\newblock In preparation. 2024.

\bibitem[NOb]{nekarda2}
S.~Nekarda and T.~Ohmoto.
\newblock Enumerative geometry of space curve projections using {T}hom polynomials.
\newblock In preparation. 2024.

\bibitem[Ohm03]{ohmotoCHI}
T.~Ohmoto.
\newblock An elementary remark on the integral with respect to {E}uler characteristics of projective hyperplane sections.
\newblock {\em Rep. Fac. Sci. Kagoshima Univ.}, 36:37--41, 2003.

\bibitem[Ohm06]{OhmotoCamb}
T.~Ohmoto.
\newblock Equivariant {C}hern classes of singular algebraic varieties with group actions.
\newblock {\em Math. Proc. of the Cambridge Phil. Soc.}, 140(01)(115), 2006.

\bibitem[Ohm16]{ohmotoSMTP}
T.~Ohmoto.
\newblock Singularities of maps and characteristic classes.
\newblock In {\em School on Real and Complex Singularities in S{\~a}o Carlos, 2012}, number~68 in Adv. Studies in Pure Math., pages 191--265, 2016.

\bibitem[Ohm24]{OhmotoSurvey}
T.~Ohmoto.
\newblock Thom polynomials for singularities of maps.
\newblock In preparation, 2024.

\bibitem[PR22]{sutipoj}
S.~Promtapan and R.~Rim\'anyi.
\newblock Characteristic classes of symmetric and skew-symmetric degeneracy loci.
\newblock In P.~Aluffi, D.~Anderson, M.~Hering, M.~Mustata, and S.~Payne, editors, {\em Facets of Algebraic Geometry. A Collection in Honor of {W}illiam {F}ulton’s 80th Birthday}, number 472 in LMS Lecture Note Series, pages 254--283. CUP, 2022.

\bibitem[PW07]{pragacz:positivity}
P.~Pragacz and A.~Weber.
\newblock Positivity of {S}chur function expansions of {T}hom polynomials.
\newblock {\em Fundamenta Mathematicae}, 195:85--95, 2007.

\bibitem[Rim01]{rrtp}
R.~Rim\'anyi.
\newblock Thom polynomials, symmetries and incidences of singularities.
\newblock {\em Inv. Math.}, 143:499--521, 2001.

\bibitem[Rim24]{primer}
R.~Rim{\'a}nyi.
\newblock Thom polynomials. {A} primer.
\newblock arXiv:2407.13883, 2024.

\bibitem[RV18]{RV}
R.~Rim{\'a}nyi and A.~Varchenko.
\newblock Equivariant {C}hern-{S}chwartz-{M}ac{P}herson classes in partial flag varieties: interpolation and formulae.
\newblock In {\em Schubert Varieties, Equivariant Cohomology and Characteristic Classes}, IMPANGA2015, pages 225--235. EMS, 2018.

\bibitem[Tho56]{thomoriginal}
R.~Thom.
\newblock Les singularit{\'e}s des applications diff{\'e}rentiables.
\newblock {\em Ann. Inst. Fourier}, 6:43--87, 1955-56.

\bibitem[TPP]{TPP}
R. Rim\'anyi: Thom {P}olynomial {P}ortal. An online registry of known Thom polynomials, https://tpp.web.unc.edu.

\end{thebibliography}

\end{document}